\documentclass[preprint,12pt]{elsarticle}
\usepackage{graphicx}
\usepackage{amsmath,amssymb,mathrsfs}
\newtheorem{thm}{Theorem}[section]

\newtheorem{lem}[thm]{Lemma}

\numberwithin{equation}{section}

\newproof{pf}{Proof}
\newproof{pot}{Proof of Theorem \ref{thm2}}

\usepackage{amssymb}

\journal{Journal de Math\'{e}matiques Pures et Appliqu\'{e}es}

\begin{document}

\begin{frontmatter}

\title{Enhanced Near-cloak by FSH Lining}

\author[a]{Hongyu Liu\corref{col}}
\ead{hongyu.liuip@gmail.com}

\author[b]{Hongpeng Sun\fnref{label3}}

\cortext[col]{Corresponding author.}

\address[a]{Department of Mathematics and Statistics, University of North Carolina, Charlotte, NC 28263, USA.}
\ead{hpsun@amss.ac.cn}
\address[b]{Institute of Mathematics, Academy of Mathematics and
Systems Science, Chinese Academy of Sciences, Beijing 100190, P. R.
China.\fnref{label3}}
\fntext[label3]{The work of this author was partially
supported by grant under NSF No.10990012.}
\begin{abstract}

We consider regularized approximate cloaking for the Helmholtz equation. Various cloaking schemes have been recently proposed and extensively investigated. The existing cloaking schemes in literature are (optimally) within
$|\ln\rho|^{-1}$ in 2D and $\rho$ in 3D of the perfect cloaking, where $\rho$ denotes the regularization
parameter. In this work, we develop a cloaking scheme with a well-designed lossy layer right outside the cloaked region that can produce significantly enhanced near-cloaking performance. In fact, it is proved that the proposed cloaking scheme could (optimally) achieve
$\rho^N$ in $\mathbb{R}^N$, $N\geq 2$, within the perfect cloaking. It is also shown that the proposed lossy layer is a finite realization of a sound-hard layer. We work with general geometry and arbitrary cloaked contents of the proposed cloaking device.

%

\medskip

\noindent{\bf R\'esum\'e}

\medskip

Nous consid\'erons le probl\`eme d'invisibilit\'e approch\'ee  pour l'\'equation d'Helmholtz.
Diverses m\'ethodes ont \'et\'e r\'ecemment propos\'ees et \'etudi\'ees.
Les techniques de quasi-invisibilit\'e pr\'esentes dans la litt\'erature approchent
l'invisibilit\'e parfaite avec une erreur proportionelle \`a $| \ln \rho |^{-1}$ dans $\mathbb{R}^2$
et $\rho$ dans $\mathbb{R}^3$, o\`u $\rho$ d\'esigne le param\`etre de r\'egularisation.
Dans ce travail, nous d\'eveloppons un syst\`eme d'invisibilit\'e qui utilise une couche
avec perte \`a l'ext\'erieur de la r\'egion dissimul\'ee et am\'eliore consid\'erablement
la quasi-invisibilit\'e.
Nous prouvons que cette nouvelle technique de dissimulation approche l'invisibilit\'e parfaite
avec une erreur proportionelle \`a $\rho^N$ dans $\mathbb{R}^N$, $N \geq 2$.
Nous prouvons \'egalement que cette couche avec perte est un cas particulier d'une couche rigide.
Ce travail concerne des dispositifs de dissimulation avec une g\'eom\'etrie g\'en\'erale.

\medskip


\end{abstract}

\end{frontmatter}


\section{Introduction}

A region is said to be \emph{cloaked} if its contents together with
the cloak are invisible to certain measurements. From a practical
viewpoint, these measurements are made in the exterior of the cloak.
Blueprints for making objects invisible to electromagnetic waves
were proposed by Pendry {\it et al.} \cite{PenSchSmi} and Leonhardt
\cite{Leo} in 2006. In the case of electrostatics, the same idea was
discussed by Greenleaf {\it et al.} \cite{GLU2} in 2003. The key
ingredient is that optical parameters have transformation properties
and could be {\it pushed-forward} to form new material parameters.
The obtained materials/media are called {\it transformation media}.
We refer to \cite{CC,GKLU4,GKLU5,Nor,U2,YYQ} for state-of-the-art
surveys on the rapidly growing literature and many striking
applications of the so-called `transformation optics'.

In this work, we shall be mainly concerned with the cloaking for the time-harmonic scalar waves governed by the
Helmholtz equation. The transformation media
proposed in \cite{GLU2,PenSchSmi} are rather singular. This poses
much challenge to both theoretical analysis and practical
fabrication. In order to avoid the singular structures, several
regularized approximate cloaking schemes are proposed in
\cite{GKLU_2,KOVW,KSVW,Liu,RYNQ}. The idea is either to incorporate
regularization into the singular transformation underlying the ideal
cloaking, or to truncate a thin layer of the singular cloaking medium near the cloaking interface.
Instead of the perfect invisibility, one would
consider the `near-invisibility' depending on a regularization
parameter. Our study is closely related to the one introduced in
\cite{KSVW} for approximate cloaking in electric impedance tomography, where the
`blow-up-a-point' transformation in \cite{GLU2,PenSchSmi} is
regularized to be the `blow-up-a-small-region' transformation. The
idea was further explored in \cite{KOVW,Liu,NgVo} for the
Helmholtz equation. In \cite{Liu}, the author imposed a homogeneous
Dirichlet boundary condition at the inner edge of the cloak and
showed that the `blow-up-a-small-region' construction gives
successful near-cloak. In \cite{KOVW}, the authors introduced a
special lossy-layer between the cloaked region and the cloaking
region, and also showed that the `blow-up-a-small-region'
construction gives successful near-cloak. For both cloaking
constructions, it was shown that the near-cloaks come, respectively,
within $1/|\ln \rho|$ in 2D and $\rho$ in 3D of the perfect
cloaking, where $\rho$ is the relative size of the small region
being blown-up for the construction and plays the role of a
regularization parameter. These estimates are also shown to be
optimal for their constructions. More subtle issues of the lossy-layer cloaking construction
developed in \cite{KOVW} were studied in \cite{NgVo}.

It is worth noting that if one lets the lossy parameter in
\cite{KOVW} go to infinity, this limit corresponds to the imposition
of a homogeneous Dirichlet boundary condition at the inner edge of
the cloak. On the other hand, the imposition of a homogeneous
Dirichlet boundary condition at the inner edge of the cloak is
equivalent to employing a sound-soft layer right outside the cloaked
region. In this sense, the lossy layer lining in \cite{KOVW} is a
finite realization of the sound-soft lining in \cite{Liu}. We would
like to emphasize that employing some special lining is necessary
for a successful near-cloaking construction, since otherwise it is shown in
\cite{KOVW} that there exists resonant inclusions which
defy any attempt to achieve near-cloak.

Though the existing cloaking constructions would yield successful near-cloaks, cloaking schemes with enhanced cloaking performances would
clearly be of great desire and significant practical importance, especially in the 2D case as can be seen from our earlier
discussion. A novel regularized cloaking scheme were developed in \cite{LiLiuSun} by making use of an FSH lining.
The FSH lining is a special lossy layer with well-designed material parameters. The study in \cite{LiLiuSun} is conducted for cloaking
device with spherical geometry and uniform cloaked contents, where the authors rely on spherical wave series representation of the
underlying wave field to derive the estimates of the cloaking performance. The newly developed cloaking scheme is shown to produce significantly enhanced
cloaking performance. In this work, we shall prove the general case with general geometry and arbitrary cloaked contents of the FSH lining construction.
For the construction, it is shown that one could achieve,
respectively, $\rho^2$ in 2D and $\rho^3$ in 3D within the perfect
cloaking. Apparently, our novel cloaking proposal with such significantly improved cloaking performances would be a very promising scheme for
constructing practical cloaking device. From our arguments in deriving these estimates, one can see that the FSH layer is a finite realization of a sound-hard
layer. Hence, the FSH layer is of completely different physical nature from the one in \cite{KOVW} which is a finite realization of a sound-soft layer. In fact, the one in \cite{KOVW} makes essential use of a large lossy parameter, whereas for our FSH layer we only require a finite lossy parameter but a large density parameter of the layer medium. 

The analysis of cloaking must specify the type of exterior
measurements. In \cite{GKLU_2,KOVW,KSVW}, the near-cloaks are
assessed in terms of boundary measurements encoded into the boundary Neumann-to-Dirichlet map or
Dirichlet-to-Neumann (DtN) map. The scattering measurement encoded into the scattering amplitude is
considered for the near-cloaks in \cite{LiLiuSun,Liu}. In the current article, we
shall assess our near-cloak construction with respect to the boundary measurements. Nonetheless, by
\cite{Nach,NSU}, it is known that
knowing the boundary DtN/NtD map amounts to knowing the scattering amplitude.

In this paper, we focus entirely on transformation-optics-approach
in constructing cloaking devices. But we would
like to mention in passing the other promising cloaking schemes
including the one based on anomalous localized resonance \cite{MN},
and another one based on special (object-dependent) coatings
\cite{AE}. {It is also interesting to note a
recent work in \cite{AKLL}, where the authors implement multi-coatings to enhance the near-cloak
in EIT. The same idea has also been extended to acoustic cloaking for achieving
enhancement in \cite{AKLL1,AGJKLH}.}

The rest of the paper is organized as follows. In Section 2, we
develop the cloaking scheme by employing the FSH lining and present the main theorems. Section 3 is devoted to the proofs of the main results. In Section 4, we
derive some crucial estimates on small inclusions that were needed in Section 3. In Section 5, we consider our cloaking construction within spherical geometry
and uniform cloaked contents, which illustrates the sharpness of our estimates in Section 3. Section 6 is devoted to discussion.

\section{Near-cloak with FSH lining}


Let $q\in L^\infty(\mathbb{R}^N)$ be a real scalar function and $\sigma=(\sigma^{ij})_{i,j=1}^N\in\mbox{Sym}(N)$ be a
symmetric-matrix-valued function on $\mathbb{R}^N$, which is bounded in
the sense that, for some constants $0<c_0<C_0<\infty$,
\begin{equation}
\label{eqn:Bound_Sigma} c_0 \xi^T \xi \leq \xi^T \sigma(x) \xi \leq
C_0 \xi^T \xi
\end{equation}
for all $x\in \mathbb{R}^N$ and $\xi \in \mathbb{R}^N$. In acoustics, $\sigma^{-1}$ and $q$, respectively, represent the mass density tensor and the bulk modulus of a {\it
regular} acoustic medium. We shall denote
$\{\mathbb{R}^N; \sigma, q\}$ an acoustic medium as described above.
It is assumed that the inhomogeneity of the medium is
compactly supported, namely, $\sigma=I$ and $q=1$ in
$\mathbb{R}^N\backslash\bar{D}$ with $D$ a bounded Lipschitz domain
in $\mathbb{R}^N$. In $\mathbb{R}^N$, the scalar wave propagation is govern by
\[
q(x)U_{tt}-\nabla\cdot(\sigma(x)\nabla U)=0\quad\mbox{in\ \ $\mathbb{R}^N$}.
\]
The time-harmonic solutions $U(x,t)=u(x)e^{-i\omega t}$ is described by the heterogeneous Helmholtz equation
\begin{equation}\label{eq:Helm}
\sum_{i,j=1}^N\frac{\partial}{\partial x_i}\left(\sigma^{ij}(x)\frac{\partial u}{\partial x_j}\right)+\omega^2 q(x) u=0\quad\mbox{in\ \ $\mathbb{R}^N$}.
\end{equation}
Let $\Omega\subset\mathbb{R}^n$ be a bounded Lipschitz domain such that $D\subseteq\Omega$. An important problem arising from practical applications is described as following. Let $\nu=(\nu_i)_{i=1}^N$ be the exterior unit normal vector to $\partial\Omega$. Impose the following boundary condition on $\partial\Omega$ for (\ref{eq:Helm}),
\begin{equation}\label{eq:N}
\sum_{i,j=1}^N\nu_i\sigma^{ij}\frac{\partial u}{\partial x_j}=\psi\in H^{-1/2}(\partial\Omega)\quad\mbox{on\ \ $\partial\Omega$},
\end{equation}
and define the Neumann-to-Dirichlet (NtD) map by
\begin{equation}\label{eq:NtD}
\Lambda_{\sigma,q}(\psi)=u|_{\partial\Omega}\in H^{1/2}(\partial\Omega),
\end{equation}
where $u\in H^1(\Omega)$ solves (\ref{eq:Helm})--(\ref{eq:N}). It is known that $\Lambda_{\sigma,q}: H^{-1/2}(\partial\Omega)\rightarrow H^{1/2}(\partial\Omega)$ is well-defined and invertible provided $\omega^2$ avoids a discrete set of eigenvalues. The practical problem is to recover $\{D;\sigma,q\}$ by knowledge of $\Lambda_{\sigma,q}$ which encodes the exterior boundary measurements.

In this paper, we shall be concerned with the construction of a layer of cloaking medium which makes the inside medium invisible to exterior measurements. To that end, we present a quick discussion on transformation acoustics. Let $\tilde
x=F(x):\Omega\rightarrow\widetilde\Omega$ be a bi-Lipschitz and
orientation-preserving mapping. For an acoustic medium
$\{\Omega;\sigma,q\}$, we let the {\it push-forwarded} medium be
defined by
\begin{equation}\label{eq:pushforward}
\{\widetilde\Omega;\widetilde\sigma,\widetilde
q\}=F_*\{\Omega;\sigma,q\}:=\{\Omega; F_*\sigma, F_*q\},
\end{equation}
where
\begin{equation}\label{tranform}
\begin{split}
&\widetilde{\sigma}(\tilde
x)=F_*\sigma(x):=\frac{1}{J}M\sigma(x)M^T|_{x=F^{-1}(\tilde x)}\\
&\widetilde{q}(\tilde x)=F_*q(x):=q(x)/J|_{x=F^{-1}(\tilde x)}
\end{split}
\end{equation}
and $M=(\partial \tilde{x}_i/\partial x_j)_{i,j=1}^N$,
$J=\mbox{det}(M)$. Then $u\in H^1(\Omega)$ solves the Helmholtz equation
\[
\nabla\cdot(\sigma(x)\nabla u)+\omega^2q(x) u=0\quad\mbox{on\ $\Omega$},
\]
if and only if the pull-back field $\widetilde u=(F^{-1})^*u:=u\circ F^{-1}\in H^1(\widetilde{\Omega})$ solves
\[
\widetilde{\nabla}\cdot(\widetilde\sigma(\tilde x)\widetilde\nabla \widetilde
u)+\omega^2\widetilde q(\tilde x)\widetilde u=0.
\]
We have made use of $\nabla$ and $\widetilde\nabla$ to distinguish the
differentiations respectively in $x$- and $\tilde x$-coordinates. We
refer to \cite{KOVW, Liu} for a proof of this invariance.

We are in a position to construct the cloaking device. In the sequel, we let $\Omega$ be a connected smooth domain and $D$ be a convex smooth domain, and suppose that $D\Subset\Omega$ and $\Omega\backslash\overline{D}$
is connected. W.L.O.G., we assume that $D$ contains the origin. Let $\rho>0$ be sufficiently small and $D_\rho:=\{\rho x; x\in D\}$. Suppose
\begin{equation}\label{eq:trans}
F_\rho: \overline{\Omega}\backslash D_\rho\rightarrow \overline{\Omega}\backslash D,
\end{equation}
which is a bi-Lipschitz and orientation-preserving mapping, and $F_\rho|_{\partial \Omega}=\mbox{Identity}$. A celebrated example of such blow-up mapping is given by
\begin{equation}\label{eq:F:ball:map}
y=F_\rho(x):=\left(\frac{R_1-\rho}{R_2-\rho}R_2+\frac{R_2-R_1}{R_2-\rho}|x|\right)\frac{x}{|x|},\
\ \rho<R_1<R_2
\end{equation}
which blows-up the central ball $B_\rho$ to $B_{R_1}$ within $B_{R_2}$. Now, we set
\begin{equation}\label{eq:F}
F(x)=\begin{cases}
\ F_\rho(x)\quad &\mbox{for\ $x\in\Omega\backslash D_\rho$},\\
\ \frac{x}{\rho} \quad &\mbox{for\ $x\in D_\rho$}.
\end{cases}
\end{equation}
Clearly, $F: \Omega\rightarrow\Omega$ is bi-Lipschitz and orientation-preserving and $F|_{\partial \Omega}=\mbox{Identity}$.
Next, let
\begin{equation}\label{eq:lossy virtual}
\{D_\rho\backslash\overline{D}_{\rho/2}; \sigma_l, q_l\},\quad \sigma_l=\gamma\rho^{2+\delta} I,\ q_l=\alpha+i\beta,
\end{equation}
where $\alpha,\beta, \gamma, \delta$ are fixed positive constants, and
\begin{equation}\label{eq:lossy physical}
\{D\backslash\overline{D}_{1/2}; \sigma_l', q_l'\}=F_*\{D_\rho\backslash\overline{D}_{\rho/2}; \sigma_l, q_l\}.
\end{equation}
We further let
\begin{equation}\label{eq:cloaking medium}
\{\Omega\backslash\overline{D};\sigma_c^\rho,q_c^\rho\}=(F_\rho)_*\{\Omega\backslash\overline{D}_{\rho};
I, 1\}.
\end{equation}
Let $D_{1/2}$ represent the region which we intend to cloak and
\begin{equation}\label{eq:target physical}
\{D_{1/2}; \sigma_a', q_a'\}
\end{equation}
be the target medium which is {\it arbitrary} but {\it regular}. We claim the following construction yields a near-cloaking device
occupying $\Omega$,
\begin{equation}\label{eq:cloaking device}
\{\Omega; \sigma, q\}=\begin{cases}
\ \ \sigma_c^\rho, q_c^\rho\quad & \mbox{in\ $\Omega\backslash\overline{D}$},\\
\ \ \sigma_l', q_l'\quad & \mbox{in\ $D\backslash\overline{D}_{1/2}$},\\
\ \ \sigma_a', q_a'\quad & \mbox{in\ $D_{1/2}$}.
\end{cases}
\end{equation}

In order to present the main theorem justifying the near-cloaking construction (\ref{eq:cloaking device}), we let $u_0$ be a solution to the following PDE system
\begin{equation}\label{eq:free space}
\begin{cases}
\Delta u_0+\omega^2 u_0=0\quad & \mbox{in\ $\Omega$},\\
\displaystyle{\frac{\partial u_0}{\partial\nu}=\psi}\quad & \mbox{on\ $\partial\Omega$}.
\end{cases}
\end{equation}
That is, $u_0$ is the wave field in the ``free space". We suppose that $-\omega^2$ is not an eigenvalue of the Neumann Laplacian. Hence, one
has a well-defined ``free" NtD map
\[
\Lambda_0(\psi)=u_0|_{\partial\Omega},
\]
where $u_0$ solves (\ref{eq:free space}). We have

\begin{thm}\label{thm:main}
Suppose $-\omega^2$ is not an eigenvalue of the Laplacian on $\Omega$ with Neumann boundary condition. Let $\Lambda_{\sigma,q}$ be the NtD map corresponding to the construction (\ref{eq:cloaking device}), and $\Lambda_0$ be the ``free" NtD map. Then there exists a constant $\rho_0$ such that for any $\rho<\rho_0$,
\begin{equation}\label{eq:main1}
\left\|\Lambda_{\sigma,q}-\Lambda_0\right\|_{\mathcal{L}(H^{-1/2}(\partial\Omega), H^{1/2}(\partial\Omega))}\leq C\rho^N,
\end{equation}
where $C$ is a positive constant dependent only on $\rho_0$, $\omega$, $\alpha$, $\beta$, $\gamma$ and $D$, $\Omega$, but completely independent of $\rho$. That is, the construction (\ref{eq:cloaking device}) produce a near-cloaking scheme which is within $\rho^N$ of the perfect cloaking in $\mathbb{R}^N$.
\end{thm}

\section{Proof of the main result}\label{sect:3}

This section is devoted to the proof of Theorem~\ref{thm:main}. First, for $\{\Omega; \sigma, q\}$ given in (\ref{eq:cloaking device}), we let
\begin{equation}\label{eq:virtual}
\{\Omega; \sigma_\rho, q_\rho\}=(F^{-1})_*\{\Omega; \sigma, q\}=\begin{cases}
\ \ I, 1\quad & \mbox{in \ $\Omega\backslash\bar{D}_\rho$},\\
\ \ \sigma_l, q_l\quad & \mbox{in \ $D_\rho\backslash\bar{D}_{\rho/2}$},\\
\ \ \sigma_a, q_a\quad & \mbox{in \ $D_{\rho/2}$},
\end{cases}
\end{equation}
where
\[
\{D_{\rho/2};\sigma_a, q_a\}=(F^{-1})_*\{D_{1/2}; \sigma_a', q_a'\}.
\]
We consider the solution $u_\rho$ of
\begin{equation}\label{eq:virtual wave}
\begin{cases}
\ \nabla\cdot(\sigma_\rho\nabla u_\rho)+\omega^2 q_\rho u_\rho=0\quad  & \mbox{in\ \ $\Omega$},\\
\ \displaystyle{\frac{\partial u_\rho}{\partial\nu}=\psi}\quad & \mbox{on\ $\partial\Omega$}.
\end{cases}
\end{equation}
Noting $F|_{\partial\Omega}=\mbox{Identity}$, by the transformation acoustics, it is straightforward to show that
\begin{equation}\label{eq:fact}
\Lambda_{\sigma, q}(\psi)=\Lambda_{\sigma_\rho, q_\rho}(\psi),\quad \forall \psi\in H^{-1/2}(\partial\Omega).
\end{equation}
Hence, in order to prove Theorem~\ref{thm:main}, we only need to show

\begin{thm}\label{thm:1}
Suppose $-\omega^2$ is not an eigenvalue of the Laplacian on $\Omega$ with Neumann boundary condition. Let $u_0$ and $u_\rho$ be the solutions of (\ref{eq:free space}) and (\ref{eq:virtual wave}) respectively. Then there exists a constant $\rho_0>0$ such that for any $\rho<\rho_0$,
\begin{equation}\label{eq:3}
\|u_\rho-u_0\|_{H^{1/2}(\partial\Omega)}\leq C\rho^N\|\psi\|_{H^{-1/2}(\partial\Omega)},
\end{equation}
where $C$ is a constant dependent only on $\rho_0$, $\omega$, $\alpha$, $\beta$, $\gamma$ and $D$, $\Omega$, but independent of $\rho$ and $\psi$.
\end{thm}

Our proof of Theorem~\ref{thm:1} would follow the spirit of the one for proving the main theorem in \cite{KOVW}. However, the main strategy in \cite{KOVW} is to control the Dirichlet value of $u_\rho$ on the exterior of the lossy layer, namely $\partial D_\rho$, and then derive some estimates of exterior boundary effects due to small sound-soft like inclusions; whereas in our case, we would control the value of the conormal derivative of $u_\rho$ on the exterior of the lossy layer $\partial D_\rho^+$, and then derive some estimates of exterior boundary effects due to small sound-hard like inclusions. It is also emphasized that by making use of layer potential techniques, we work with general geometry of the cloaking device.

We first derive the following lemma.

\begin{lem}\label{lem:1}
The solutions of (\ref{eq:free space}) and (\ref{eq:virtual wave}) satisfy
\begin{equation}\label{eq:4}
\beta\omega^2\int_{D_{\rho}\backslash\overline{D}_{\rho/2}}|u_\rho|^2\ dx\leq C\|\psi\|_{H^{-1/2}(\partial\Omega)}\|u_\rho-u_0\|_{H^{1/2}(\partial\Omega)},
\end{equation}
where $C$ is a positive constant (depending only on $\Omega$).
\end{lem}
\begin{pf}
Multiplying (\ref{eq:virtual wave}) by $\bar{u}_\rho$ and integrating by parts, we have
\begin{equation}\label{eq:5}
-\int_{\Omega} \sigma_\rho |\nabla u_\rho|^2\ dx+\omega^2\int_{\Omega} q_\rho |u_\rho|^2\ dx=-\int_{\partial\Omega} (\sigma_\rho\nabla u_\rho)\cdot \nu \bar{u}_\rho\ d\sigma_x,
\end{equation}
which in turn yields
\begin{equation}\label{eq:6}
\begin{split}
&\beta\omega^2\int_{D_{2\rho}\backslash\overline{D}_\rho}|u_\rho|^2\ dx\\
=&-\Im\left(\int_{\partial\Omega}\frac{\partial u_\rho}{\partial\nu}\cdot \bar{u}_\rho\ d\sigma_x\right)=-\Im\left(\int_{\partial\Omega}\psi(\bar{u}_\rho-\bar{u}_0)\ d\sigma_x\right).
\end{split}
\end{equation}
By (\ref{eq:6}), we immediately have (\ref{eq:4}).
\end{pf}
In the following, we let
\begin{equation}\label{eq:nei}
\Psi^-(x)=\nu\cdot\nabla u_\rho^- (x)\quad \mbox{on\ $\partial D_\rho$},
\end{equation}
namely, the normal derivative of $u_\rho$ on $\partial D_{\rho}$ when one approaches $\partial D_{\rho}$ from the interior of $D_{\rho}$. Here and throughout the rest of this paper, $\nu$ denotes the exterior unit normal of the domain under discussion.
Similarly, we let
\begin{equation}\label{eq:wai}
\Psi^+(x)=\nu\cdot\nabla u_\rho^+(x)\quad \mbox{on\ $\partial D_\rho$}
\end{equation}
denote the normal derivative of $u_\rho$ on $\partial D_\rho$ when one approaches $\partial D_\rho$ from the exterior of $D_\rho$. We shall show

\begin{lem}\label{lem:m}
The solutions to (\ref{eq:free space}) and (\ref{eq:virtual wave}) verify
\begin{equation}\label{eq:imp1}
\begin{split}
&\left\|\Psi^-(\rho\ \cdot)\right\|_{H^{-3/2}(\partial D)}^2\\
\leq &\ C \frac{(\gamma+\sqrt{\alpha^2+\beta^2}\rho^{-\delta}\omega^2)^2}{\beta\gamma^2\omega^2}\rho^{-N-2}\|\psi\|_{H^{-1/2}(\partial\Omega)}\|u_\rho
-u_0\|_{H^{1/2}(\partial\Omega)},
\end{split}
\end{equation}
and
\begin{equation}\label{eq:imp11}
\begin{split}
&{\left\|\Psi^+(\rho\ \cdot)\right\|_{H^{-3/2}(\partial D)}^2}\\
\leq &\ C \frac{(\gamma+\sqrt{\alpha^2+\beta^2}\rho^{-\delta}\omega^2)^2}{\beta\omega^2}\rho^{2(1+\delta)-N}\|\psi\|_{H^{-1/2}(\partial\Omega)}\|u_\rho
-u_0\|_{H^{1/2}(\partial\Omega)},
\end{split}
\end{equation}
where $C$ is positive constant dependent only on $D$ and $\Omega$, but independent of $\psi$ and $\rho$.
\end{lem}

\begin{pf}
We shall make use of the following fact
\begin{equation}\label{eq:p1}
\left\|\Psi(\rho\ \cdot)\right\|_{H^{-3/2}(\partial D)}=\sup_{\|\phi\|_{H^{3/2}(\partial D)}\leq 1}\left|\int_{\partial D}\Psi(\rho x)\phi(x)\ d\sigma_x\right|.
\end{equation}
For any $\phi\in H^{3/2}(\partial D)$, there exists $w\in H^2(D)$ such that
\begin{align*}
& (i)~~\mbox{$w=\phi$\ on $\partial D$\ and \ $\frac{\partial w}{\partial\nu}=0$\ on $\partial D$},\\
& (ii)~~\mbox{$\|w\|_{H^2(D)}\leq C \|\phi\|_{H^{3/2}(\partial D)}$},\\
& (iii)~~\mbox{$w=0$\ in $D_{1/2}$}.
\end{align*}
Then we have
\begin{equation}\label{eq:bc}
\int_{\partial D}\Psi^-(\rho x)\phi(x)\ d\sigma_x=\int_{\partial D} \frac{\partial u_\rho^-}{\partial \nu}(\rho x)\phi(x)\ d\sigma_x=\int_{\partial D}\frac{\partial u_\rho^-}{\partial \nu}(\rho x) w(x)\ d\sigma_x.
\end{equation}
For $y\in D_{\rho}$, let
\[
x:=\frac{y}{\rho}\in D.
\]
Set
\[
v(x):=u_\rho(\rho x)=u_\rho(y),\ \ x\in D.
\]
Since
\begin{equation}\label{eq:nei1}
\gamma\nabla_y\cdot(\rho^{2+\delta}\nabla_y u_\rho)+\omega^2(\alpha+i\beta) u_\rho=0\quad\mbox{in\ \ $D_{\rho}\backslash\overline{D}_{\rho/2}$},
\end{equation}
it is directly verified that
\begin{equation}\label{eq:n1}
\gamma\nabla_x\cdot(\rho^\delta\nabla_x v)+\omega^2(\alpha+i\beta)v=0\quad\mbox{in\ \ $D\backslash\overline{D}_{1/2}$}.
\end{equation}
By Green's formula and (\ref{eq:bc}), we have
\begin{equation}\label{eq:n2}
\begin{split}
& \int_{\partial D}\Psi^-(\rho x)\phi(x)\ d\sigma_x\\
=&\int_{\partial D}\frac{\partial u_\rho^-}{\partial\nu}(\rho x)\phi(x)\ d\sigma_x
=\rho^{-1}\int_{\partial D}\frac{\partial v^-}{\partial \nu}(x)\phi(x)\ d\sigma_x\\
=& \rho^{-1}\left[\int_{\partial D} \frac{\partial v^-}{\partial \nu}(x)w(x)\ d\sigma_x-\int_{\partial D} v(x)\frac{\partial w}{\partial\nu}(x)\ d\sigma_x\right]\\
=& \rho^{-1}\left[\int_{D}\Delta v(x)w(x)\ dx-\int_{D} v(x)\Delta w(x)\ dx\right].
\end{split}
\end{equation}
Then by (\ref{eq:n1}) and (\ref{eq:n2}), we further have
\begin{equation}\label{eq:n3}
\begin{split}
& \left|\int_{\partial D}\Psi^-(\rho x)\phi(x)\ d\sigma_x\right|\\
\leq & \rho^{-1}\left|\int_{D}\Delta v(x)w(x)\ dx-\int_{D} v(x)\Delta w(x)\ dx\right|\\
\leq & \frac{\sqrt{\alpha^2+\beta^2}}{\gamma}\rho^{-\delta-1}\omega^2\left(\int_{D\backslash\overline{D}_{1/2}}|v(x)|^2 dx\right)^{1/2}\|w\|_{L^2(D)}\\
&+\rho^{-1}\left(\int_{D\backslash\overline{D}_{1/2}}|v(x)|^2 dx\right)^{1/2}\|\Delta w\|_{L^2(D)}.
\end{split}
\end{equation}
Using the relation
\[
\|v\|_{L^2(D\backslash\overline{D}_{1/2})}=\|u_\rho(\rho\ \cdot)\|_{L^2(D\backslash\overline{D}_{1/2})}=\rho^{-N/2}\|u_\rho\|_{L^2(D_{\rho}\backslash\overline{D}_{\rho/2})}
\]
we have from (\ref{eq:n3}) that
\begin{equation}\label{eq:n4}
\begin{split}
& \left|\int_{\partial D}\Psi^-(\rho x)\phi(x)\ d\sigma_x\right|\\
\leq &\ C \rho^{-N/2-1}\left(1+\frac{\sqrt{\alpha^2+\beta^2}}{\gamma}\rho^{-\delta}\omega^2\right)
\|u_\rho\|_{L^2(D_{\rho}\backslash\overline{D}_{\rho/2})}
\|\phi\|_{H^{3/2}(\partial D)},
\end{split}
\end{equation}
which implies
\begin{equation}\label{eq:n5}
\left\|\Psi^-(\rho\ \cdot)\right\|_{H^{-3/2}(\partial D)}\leq C\rho^{-N/2-1}\left(1+\frac{\sqrt{\alpha^2+\beta^2}}{\gamma}\rho^{-\delta}\omega^2\right)\|u_\rho\|_{L^2(D_{\rho}\backslash\overline{D}_{\rho/2})}.
\end{equation}
By (\ref{eq:n5}) and Lemma~\ref{lem:1}, one immediately has (\ref{eq:imp1}). Finally, by (\ref{eq:nei1}) and the transmission condition on $\partial D$, we see
\[
\frac{\partial u^+}{\partial\nu}\bigg|_{\partial D_\rho}=\gamma\rho^{2+\delta}\frac{\partial u^-}{\partial \nu}\bigg|_{\partial D_\rho}
\]
and hence
\[
\Psi^+(\rho x)=\gamma\rho^{2+\delta} \Psi^-(\rho x)\quad \mbox{for \ $x\in \partial D$}
\]
which together with (\ref{eq:imp1}) implies (\ref{eq:imp11}).

The proof is completed.

\end{pf}

The next lemma is of crucial importance in proving Theorem~\ref{thm:1}.

\begin{lem}\label{lem:uniform wellposed}
Suppose $-\omega^2$ is not an eigenvalue of the Laplacian on $\Omega$ with Neumann boundary condition. Let $u_0\in H^1(\Omega)$ be the solution of (\ref{eq:free space}). Let $\varphi\in H^{-1/2}(\partial D_\tau)$ and consider the Helmholtz system
\begin{equation}\label{eq:uni 1}
\begin{cases}
\ \ \displaystyle{\Delta u_\tau+\omega^2 u_\tau=0}\quad & \mbox{in\ \ $\Omega\backslash \overline{D}_\tau$},\\
\ \ \displaystyle{\frac{\partial u_\tau}{\partial \nu}=\varphi}\quad & \mbox{on\ \ $\partial D_\tau$},\\
\ \ \displaystyle{\frac{\partial u_\tau}{\partial \nu}=\psi}\quad & \mbox{on\ \ $\partial\Omega$}.
\end{cases}
\end{equation}
Let
\[
\varphi_0(x)=\frac{\partial u_0}{\partial\nu}(x)\quad\mbox{for\ $x\in\partial D_\tau$}.
\]
Then there exist a constant $\tau_0>0$ such that for any $\tau<\tau_0$,
\begin{equation}\label{eq:imp3}
\|u_\tau-u_0\|_{H^{1/2}(\partial\Omega)}\leq C\ \left(\tau^N\|\psi\|_{H^{-1/2}(\partial \Omega)}+\tau^{N-1}\ \|\varphi(\tau\ \cdot)\|_{H^{-3/2}(\partial D)}\right),
\end{equation}
where $C$ is a positive constant dependent only on $\tau_0$, $\omega$ and $\Omega$, $D$, but independent of $\tau$ and $\varphi$, $\psi$.
\end{lem}

\begin{pf}
Let
\begin{equation}\label{eq:def1}
V=u_\tau-u_0\quad \mbox{on\ \ $\Omega\backslash\overline{D}_\tau$}.
\end{equation}
By (\ref{eq:free space}) and (\ref{eq:uni 1}), one sees that $V\in H^1(\Omega\backslash\overline{D}_\tau)$ satisfies
\begin{equation}\label{eq:V}
\begin{cases}
\ \ \Delta V+\omega^2 V=0\quad & \mbox{in\ $\Omega\backslash\overline{D}_\tau$},\\
\ \ \displaystyle{\frac{\partial V}{\partial \nu}=\varphi-\frac{\partial u_0}{\partial \nu}}\quad & \mbox{on\ $\partial D_\tau$},\\
\ \ \displaystyle{\frac{\partial V}{\partial \nu}=0}\quad & \mbox{on\ $\partial \Omega$}.
\end{cases}
\end{equation}
Let $V=V_1-V_2$ with $V_1\in H^1(\Omega\backslash\overline{D}_\tau)$ satisfying
\begin{equation}\label{eq:V1}
\begin{cases}
\ \ \Delta V_1+\omega^2 V_1=0\quad & \mbox{in\ $\Omega\backslash\overline{D}_\tau$},\\
\ \ \displaystyle{\frac{\partial V_1}{\partial \nu}=\varphi}\quad & \mbox{on\ $\partial D_\tau$},\\
\ \ \displaystyle{\frac{\partial V_1}{\partial \nu}=0}\quad & \mbox{on\ $\partial \Omega$}.
\end{cases}
\end{equation}
and $V_2\in H^1(\Omega\backslash\overline{D}_\tau)$ satisfying
\begin{equation}\label{eq:V2}
\begin{cases}
\ \ \Delta V_2+\omega^2 V_2=0\quad & \mbox{in\ $\Omega\backslash\overline{D}_\tau$},\\
\ \ \displaystyle{\frac{\partial V_2}{\partial \nu}=\varphi_0}\quad & \mbox{on\ $\partial D_\tau$},\\
\ \ \displaystyle{\frac{\partial V_2}{\partial \nu}=0}\quad & \mbox{on\ $\partial \Omega$}.
\end{cases}
\end{equation}
By Lemma~\ref{lem:order n} in Section~\ref{section:4}, we know
\begin{equation}\label{eq:order n}
\|V_2\|_{H^{1/2}(\partial\Omega)}\leq C \tau^N\|\psi\|_{H^{-1/2}(\partial \Omega)}.
\end{equation}
In order to estimate $\|V_1\|_{H^{1/2}(\partial\Omega)}$, we let $W\in H_{loc}^1(\mathbb{R}^N\backslash\overline{D}_\tau)$ be the unique solution to the following scattering problem
\begin{equation}\label{eq:s1}
\begin{cases}
\ \ \Delta W+\omega^2 W=0\quad & \mbox{in\ $\mathbb{R}^N\backslash\overline{D}_\tau$},\\
\ \  \displaystyle{\frac{\partial W}{\partial \nu}=\varphi}\quad & \mbox{on\ $\partial D_\tau$},\\
\ \  \displaystyle{\lim_{|x|\rightarrow +\infty}|x|^{(N-1)/2}\left\{\frac{\partial W}{\partial |x|}-i\omega W\right\}=0.}
\end{cases}
\end{equation}
Let $\tau_0$ be sufficiently small such that
\begin{equation}\label{eq:r0}
D_{\tau_0}\Subset B_{r_0}\Subset D
\end{equation}
for some finite $r_0>0$, where $B_r$ denotes a central ball of radius $r$. Let $r_0<r_1<r_2<+\infty$ be such that
\begin{equation}\label{eq:r12}
B_{r_1}\Subset\Omega\quad \mbox{and}\quad \Omega\backslash \overline{D}\Subset B_{r_2}\backslash \overline{B}_{r_0}.
\end{equation}
By Lemma~\ref{lem:imp1} in Section~\ref{section:4}, we have
\begin{equation}\label{eq:lemimp1}
\|W\|_{L^2(B_{r_2}\backslash\overline{B}_{r_0})}\leq C \tau^{N-1}\|\varphi(\tau\ \cdot)\|_{H^{-3/2}(\partial D)}.
\end{equation}
Since $(\Delta+\omega^2)W=0$, by the interior regularity estimates, we see
\begin{align}
&\left\|\frac{\partial W}{\partial \nu}\bigg|_{\partial B_{r_1}}\right\|_{H^{1/2}(\partial B_{r_1})}\leq C \tau^{N-1}\|\varphi(\tau\ \cdot)\|_{H^{-3/2}(\partial D)},\label{eq:lemimp2}\\
&\left\|W\right\|_{H^{1/2}(\partial B_{r_1})}\leq C \tau^{N-1}\|\varphi(\tau\ \cdot)\|_{H^{-3/2}(\partial D)},\label{eq:lemimp3}
\end{align}
and
\begin{equation}\label{eq:lemimp4}
\left\|W\right\|_{H^{1/2}(\partial \Omega)}\leq C \tau^{N-1}\|\varphi(\tau\ \cdot)\|_{H^{-3/2}(\partial D)}.
\end{equation}
Next, by the Green's representation, we know
\begin{equation}\label{eq:W}
W(x)=\int_{\partial B_{r_1}} \frac{\partial G(x-y)}{\partial\nu(y)}W(y) -G(x-y)\frac{\partial W(y)}{\partial \nu(y)}\ d\sigma_y,\quad x\in\mathbb{R}^N\backslash\overline{B}_{r_1}
\end{equation}
where
\begin{equation}\label{eq:fun}
G(x)=\frac{i}{4}\left(\frac{\omega}{2\pi|x|}\right)^{(N-2)/2}H_{(N-2)/2}^{(1)}(\omega|x|)
\end{equation}
is the outgoing Green's function. By (\ref{eq:lemimp2}), (\ref{eq:lemimp3}) and (\ref{eq:W}), it is readily seen that
\begin{equation}\label{eq:W2}
\left\|\frac{\partial W}{\partial \nu}\bigg|_{\partial\Omega}\right\|_{C(\partial \Omega)}\leq C \tau^{N-1}\|\varphi(\tau\ \cdot)\|_{H^{-3/2}(\partial D)}.
\end{equation}

Let
\[
P=V_1-W.
\]
By (\ref{eq:V}) and (\ref{eq:s1}), one sees that $P\in H^1(\Omega\backslash\overline{D}_\tau)$ satisfies
\begin{equation}\label{eq:P}
\begin{cases}
\ \ \Delta P+\omega^2 P=0\quad & \mbox{in\ $\Omega\backslash\overline{D}_\tau$},\\
\ \ \displaystyle{\frac{\partial P}{\partial \nu}=0}\quad & \mbox{on\ $\partial D_\tau$},\\
\ \ \displaystyle{\frac{\partial P}{\partial \nu}=-\frac{\partial W}{\partial \nu}}\quad & \mbox{on\ $\partial \Omega$}.
\end{cases}
\end{equation}
By Lemma~\ref{lem:small inclusion} in the following, we have
\begin{equation}\label{eq:small inclusion 1}
\left\|P\right\|_{H^{1/2}(\partial\Omega)}\leq C\left\|\frac{\partial W}{\partial \nu}\bigg|_{\partial\Omega}\right\|_{C(\partial\Omega)},
\end{equation}
which together with (\ref{eq:W2}) implies
\begin{equation}\label{eq:small inclusion 2}
\left\|P\right\|_{H^{1/2}(\partial\Omega)} \leq C \tau^{N-1}\|\varphi(\tau\ \cdot)\|_{H^{-3/2}(\partial D)}.
\end{equation}

Since $V_1=P+W$, (\ref{eq:lemimp4}) and (\ref{eq:small inclusion 2}) immediately yields that
\[
\|V_1\|_{H^{1/2}(\partial\Omega)}\leq C \tau^{N-1} \|\varphi(\tau\ \cdot)\|_{H^{-3/2}(\partial D)},
\]
which together with (\ref{eq:order n}) implies (\ref{eq:imp3}).

The proof is completed.
\end{pf}

We are in a position to present the proof of Theorem~\ref{thm:1}.

\begin{pf}[Proof of Theorem~\ref{thm:1}]
By taking $\tau=\rho$ and $\varphi=\frac{\partial u_\rho^+}{\partial\nu}|_{\partial D_\rho}$ in Lemma~\ref{lem:uniform wellposed}, we have
\begin{equation}\label{eq:thm11}
\begin{split}
& \|u_\rho-u_0\|_{H^{1/2}(\partial\Omega)}\\
\leq & C_1\left(\rho^N\left\|\psi\right\|_{H^{-1/2}(\partial \Omega)}+\rho^{N-1}\left\|\left(\frac{\partial u_\rho^+}{\partial \nu}\right)(\rho\ \cdot)\right\|_{H^{-3/2}(\partial D)}\right).
\end{split}
\end{equation}
Next, by (\ref{eq:imp11}) in Lemma~\ref{lem:m}, we have for $\epsilon>0$
\begin{equation}\label{eq:m1}
\begin{split}
\left\|\left(\frac{\partial u_\rho^+}{\partial\nu}\right)(\rho\ \cdot)\right\|_{H^{-3/2}(\partial D)} \leq & C_2 \rho^{(2-N)/2}\|\psi\|_{H^{-1/2}(\partial\Omega)}^{1/2}\|u_\rho
-u_0\|_{H^{1/2}(\partial\Omega)}^{1/2}\\
\leq & C_2 \rho^{(2-N)/2}\bigg(\frac{\rho^{(2-N)/2}\cdot\rho^{N-1}}{4\epsilon}\|\psi\|_{H^{-1/2}(\partial\Omega)}\\
& +\frac{\epsilon}{\rho^{(2-N)/2}\cdot\rho^{N-1}}
\|u_\rho-u_0\|_{H^{1/2}(\partial\Omega)}\bigg)
\end{split}
\end{equation}
From (\ref{eq:thm11}) and (\ref{eq:m1}), we further have
\begin{equation}\label{eq:m2}
\begin{split}
\|u_\rho-u_0\|_{H^{1/2}(\partial\Omega)}\leq & C_1\rho^N\|\psi\|_{H^{-1/2}(\partial\Omega)}\\
&+\frac{1}{4}C_1C_2\rho^{N-1}\frac{\rho}{\epsilon}\|\psi\|_{H^{-1/2}(\partial\Omega)}
+C_1C_2\epsilon\|u_\rho-u_0\|_{H^{1/2}(\partial\Omega)}.
\end{split}
\end{equation}
By choosing $\epsilon$ such that $C_1C_2\epsilon<1$, we see immediately from (\ref{eq:m2}) that
\[
\|u_\rho-u_0\|_{H^{1/2}(\partial\Omega)}\leq C\rho^N\|\psi\|_{H^{-1/2}(\partial\Omega)},
\]
which completes the proof.
\end{pf}

\section{Some estimates on small inclusions}\label{section:4}

In this section, we shall derive those lemmas that were needed in the proof of Lemma~\ref{lem:uniform wellposed} on the wave estimates due to small inclusions. {We would like to mention that there are a lot of results on this subject in different settings in literature, see e.g., \cite{Amm3,AK1,AK2,AVV}. We shall derive some new estimates in the specific setting of our current study. We would make essential use of the layer
potential techniques to derive the desired estimates in this section}. To that end, we let $G(x)$ be the outgoing Green's function in (\ref{eq:fun}). It is well-known that when $N=2$,
\begin{equation}\label{eq:2D G}
G(x)=-\frac{1}{2\pi}\ln|x|+\frac{i}{4}-\frac{1}{2\pi}\ln\frac{\omega}{2}-\frac{E}{2\pi}+\mathcal{O}(|x|^2\ln|x|)
\end{equation}
for $|x|\rightarrow 0$, where $E$ is the Euler's constant; and when $N=3$
\begin{equation}\label{eq:3D G}
G(x)=\frac{e^{i\omega|x|}}{4\pi |x|}.
\end{equation}
For surface densities $\psi(x)$ with $x\in\partial\Omega$, and $\varphi(x)$ with $x\in\partial D_\tau$, we introduce the single- and double-layer potential operators as follows
\begin{equation}\label{eq:single}
\begin{split}
&(SL_{[\partial\Omega]}\psi)(x)=\int_{\partial\Omega} G(x-y)\psi(y)\ d\sigma_y,\ \ x\in\mathbb{R}^N\backslash\partial\Omega\\
&(SL_{[\partial D_\tau]}\varphi)(x)=\int_{\partial D_\tau} G(x-y)\varphi(y)\ d\sigma_y,\ \ x\in\mathbb{R}^N\backslash\partial D_\tau
\end{split}
\end{equation}
and
\begin{equation}\label{eq:double}
\begin{split}
&(DL_{[\partial\Omega]}\psi)(x)=\int_{\partial\Omega} \frac{\partial G(x-y)}{\partial\nu(y)}\psi(y)\ d\sigma_y,\ \ x\in\mathbb{R}^N\backslash\partial\Omega\\
&(DL_{[\partial D_\tau]}\varphi)(x)=\int_{\partial D_\tau} \frac{\partial G(x-y)}{\partial\nu(y)}\varphi(y)\ d\sigma_y,\ \ x\in\mathbb{R}^N\backslash\partial D_\tau.
\end{split}
\end{equation}
We also let
\begin{equation}\label{eq:op s}
\begin{split}
&(S_{[\partial\Omega]}\psi)(x)=\int_{\partial\Omega} G(x-y)\psi(y)\ d\sigma_y,\ \ x\in\partial\Omega\\
& (S_{[\partial D_\tau]}\varphi)(x)=\int_{\partial D_\tau} G(x-y)\varphi(y)\ d\sigma_y,\ \ x\in\partial D_\tau
\end{split}
\end{equation}
and
\begin{equation}\label{eq:op k}
\begin{split}
&(K_{[\partial\Omega]}\psi)(x)=\int_{\partial\Omega} \frac{\partial G(x-y)}{\partial\nu(y)}\psi(y)\ d\sigma_y,\ \ x\in\partial\Omega\\
& (K_{[\partial D_\tau]}\varphi)(x)=\int_{\partial D_\tau} \frac{\partial G(x-y)}{\partial\nu(y)}\varphi(y)\ d\sigma_y,\ \ x\in\partial D_\tau.
\end{split}
\end{equation}

As specified earlier in Section 3, we would like to emphasize again here that in the above integral operators, $\nu$ denotes the exterior unit normal vector of the underlying domain.

\begin{lem}\label{lem:order n}
Suppose $-\omega^2$ is not an eigenvalue of the Laplacian on $\Omega$ with Neumann boundary condition. Let $u_0\in H^1(\Omega)$ be the solution of (\ref{eq:free space}) and
$\varphi_0(x)=\frac{\partial u_0}{\partial\nu}(x)\ \mbox{for\ $x\in\partial D_\tau$}.$ Consider the Helmholtz system
\begin{equation}\label{eq:small1}
\begin{cases}
\ \ \displaystyle{\Delta u_\tau+\omega^2 u_\tau=0}\quad & \mbox{in\ \ $\Omega\backslash \overline{D}_\tau$},\\
\ \ \displaystyle{\frac{\partial u_\tau}{\partial \nu}=\varphi_0}\quad & \mbox{on\ \ $\partial D_\tau$},\\
\ \ \displaystyle{\frac{\partial u_\tau}{\partial \nu}=0}\quad & \mbox{on\ \ $\partial\Omega$}.
\end{cases}
\end{equation}
Then there exists a constant $\tau_0>0$ such that for any $\tau<\tau_0$, (\ref{eq:small1}) has a unique solution $u_\tau\in H^1(\Omega\backslash\overline{D}_\tau)$ and moreover
\begin{equation}\label{eq:small2}
\|u_\tau\|_{H^{1/2}(\partial\Omega)}\leq C\ \tau^N\|\psi\|_{H^{-1/2}(\partial \Omega)}.
\end{equation}
where $C$ is a positive constant dependent only on $\tau_0$, $\omega$ and $\Omega$, $D$, but independent of $\tau$ and $\psi$.
\end{lem}

\begin{pf}
Since $\varphi_0\in C(\partial D_\tau)$, we know $u_\tau \in C^2(\Omega\backslash\overline{D}_\tau)\cap C(\overline{\Omega}\backslash D_\tau)$ is a strong solution (cf. \cite{ColKre}). By Green's representation formula, we know
\begin{equation}\label{eq:rep1}
\begin{split}
u_\tau(x)=\int_{\partial(\Omega\backslash\overline{D}_\tau)}\left\{G(x-y)\frac{\partial u_\tau(y)}{\partial\nu(y)}-\frac{\partial G(x-y)}{\partial\nu(y)}u_\tau(y)\right\}\ d \sigma_y,\ x\in\Omega\backslash\overline{D}_\tau.
\end{split}
\end{equation}
Let
\[
h(x):=-\int_{\partial D_\tau} G(x-y)\varphi_0(y)\ d\sigma_y
\]
and
\[
\phi_1=u_\tau|_{\partial D_\tau}\quad\mbox{and}\quad \phi_2=u_\tau|_{\partial \Omega}.
\]
From (\ref{eq:rep1}) we have
\begin{equation}\label{eq:rep2}
u_\tau(x)=(DL_{[\partial D_\tau]}\phi_1)(x)-(DL_{[\partial\Omega]}\phi_2)(x)+h(x),\ \ x\in\Omega\backslash\overline{D}_\tau.
\end{equation}
Letting $x$ go to $\partial \Omega$ and $\partial D_\tau$, respectively, by the mapping properties of double-layer potential operator (cf. \cite{ColKre} and \cite{McL}), we have from (\ref{eq:rep2}) the following system of integral equations for $\phi_1\in C(\partial D_\tau)$ and $\phi_2\in C(\partial\Omega)$,
\begin{equation}\label{eq:system1}
\begin{cases}
\displaystyle{\frac{1}{2}\phi_1(x)}=&\displaystyle{(K_{[\partial D_\tau]}\phi_1)(x)-(DL_{[\partial\Omega]}\phi_2)(x)+h(x),\quad x\in\partial D_\tau}\\
\displaystyle{\frac{1}{2}\phi_2(x)}=&\displaystyle{(DL_{[\partial D_\tau]}\phi_1)(x)-(K_{[\partial\Omega]}\phi_2)(x)+h(x),\quad x\in\partial \Omega}.
\end{cases}
\end{equation}

Next, we claim
\begin{equation}\label{eq:small3}
\|h(\tau\ \cdot)\|_{C(\partial D)}\leq C\tau\|\psi\|_{H^{-1/2}(\partial\Omega)}\ \mbox{and}\ \|h(\cdot)\|_{C(\partial \Omega)}\leq C\tau^{N}\|\psi\|_{H^{-1/2}(\partial\Omega)}.
\end{equation}
We first estimate $h(x)$ for $x\in\partial\Omega$. It is noted that $|G(x-y)|,\ |\nabla_y G(x-y)|\leq \widetilde{C}$ for any $x\in\partial\Omega$ and $y\in\partial D_\tau$, where $\widetilde{C}$ depends only on $\tau_0$ and $\partial\Omega$. Integrating by parts, we have for any $x\in\partial\Omega$
\begin{equation}\label{eq:small5n}
\begin{split}
|h(x)|=&\left|\int_{\partial D_\tau}G(x-y)\frac{\partial u_0(y)}{\partial\nu(y)}\ d\sigma_y\right|\\
=& \left|\int_{D_\tau}\Delta_y u_0(y) G(x-y)\ dy+\int_{D_\tau}\nabla_yG(x-y)\cdot\nabla_y u_0(y)\ dy\right|\\
\leq & \ \omega^2\left|\int_{D_\tau}G(x-y)u_0(y)\ dy\right|+\left|\int_{D_\tau}\nabla_yG(x-y)\cdot\nabla_y u_0(y)\ dy\right|\\
\leq & \widetilde{C}\tau^N\left(\omega^2\|u_0\|_{L^\infty(D_\tau)}+\|\nabla u_0\|_{L^\infty(D_\tau)}\right)\\
\leq & C\tau^N\|\psi\|_{H^{-1/2}(\partial\Omega)}.
\end{split}
\end{equation}
We proceed to estimate $h(x)$ for $x\in\partial D_\tau$. By taking $\tau_0$ sufficiently small, we assume $D_{2\tau_0}\Subset\Omega$. We first let $x\in\partial D_{\delta\tau}$ with $1<\delta<2$. Similar to (\ref{eq:small5n}), by integration by parts, we have
\begin{equation}\label{eq:small5}
|h(x)|\leq\ \omega^2\left|\int_{D_\tau}G(x-y)u_0(y)\ dy\right|+\left|\int_{D_\tau}\nabla_yG(x-y)\cdot\nabla_y u_0(y)\ dy\right|.
\end{equation}
For the first term in (\ref{eq:small5}), we have for $x'\in \partial D_\delta$
\begin{equation}\label{eq:small50}
\begin{split}
|h_1(\tau x')|& \leq \omega^2\int_{D}\left|G(\tau(x'-y'))u_0(\tau y') \right| \tau^N dy'
\end{split}
\end{equation}
Using (\ref{eq:2D G}) and (\ref{eq:3D G}), and the mapping property of volume potential operator, it is straightforward to verify that
\begin{equation}\label{eq:small51}
\|h_1(\tau\ \cdot)\|_{C(\partial D_{\delta})}\leq \widetilde{C}\tau\|u_0(\tau\ \cdot)\|_{L^\infty(D)}\leq C\tau\|\psi\|_{H^{-1/2}(\partial\Omega)},
\end{equation}
where $C$ is independent of $\delta$ and $\tau$.
In like manner, it can be shown that for the second term in (\ref{eq:small5})
\[
h_2(x)=\int_{D_\tau}\nabla_yG(x-y)\cdot\nabla_y u_0(y)\ dy,
\]
we have
\begin{equation}\label{eq:small52}
\|h_2(\tau\ \cdot)\|_{C(\partial D_\delta)}\leq C\tau\|\psi\|_{H^{-1/2}(\partial \Omega)}.
\end{equation}
By the mapping properties of single-layer potential operator (cf. \cite{ColKre}), we know
\[
h(x)|_{\partial D_\tau}=\lim_{\delta\rightarrow 1^+}\left(h(x)|_{\partial D_{\delta \tau}}\right),
\]
which together with (\ref{eq:small51}) and (\ref{eq:small52}) implies
\begin{equation}\label{eq:small53}
\|h(\tau\ \cdot)\|_{C(\partial D)}\leq C\tau\|\psi\|_{H^{-1/2}(\partial\Omega)}.
\end{equation}

We now consider the system of integral equations (\ref{eq:system1}). It is first noted that
\begin{equation}\label{eq:small60}
\|\left(DL_{[\partial D_\tau]}\phi_1\right)(\cdot)\|_{C(\partial \Omega)}\leq C\tau^{(N-1)/2}\|\phi_1(\cdot)\|_{L^2(\partial D_\tau)}.
\end{equation}
It is further noted that since $-\omega^2$ is not an eigenvalue of the Laplacian on $\Omega$ with Neumann boundary condition,
\[
\frac{1}{2}I+K_{[\partial\Omega]}: C(\partial\Omega)\rightarrow C(\partial\Omega)
\]
is invertible (cf. \cite{ColKre}). Hence, by the second equation in (\ref{eq:system1}), and (\ref{eq:small3}), (\ref{eq:small60}), we have
\begin{equation}\label{eq:small61}
\|\phi_2\|_{C(\partial\Omega)}\leq C\left(\tau^{(N-1)/2}\|\phi_1(\cdot)\|_{L^2(\partial D_\tau)}+\tau^N\|\psi\|_{H^{-1/2}(\partial\Omega)}\right).
\end{equation}
We proceed to treat the first equation in (\ref{eq:system1}). First, by change of variables in integrals, it is straightforward to show that
\begin{equation}\label{eq:small70}
\left(K_{[\partial D_\tau]}\phi_1\right)(\tau x')=\left({K_0}_{[\partial D]}\phi_1(\tau\ \cdot)\right)(\tau x')+\left(\mathscr{R}\phi_1(\tau\ \cdot)\right)(\tau x'),\quad x'\in\partial D,
\end{equation}
where ${K_0}_{[\partial D]}$ is an integral operator with the kernel given by
${\partial G_0(x'-y')}/{\partial \nu(y')}$ with
\begin{equation}\label{eq:g0}
G_0(x'-y')=\begin{cases}
\displaystyle{-\frac{1}{2\pi}\ln|x'-y'|}\quad & N=2,\\
\displaystyle{\frac{1}{4\pi}\frac{1}{|x'-y'|}}\quad & N=3,
\end{cases}
\end{equation}
for $x'\neq y'$ and $x', y'\in\partial D$,
and $\mathscr{R}$ satisfies
\begin{equation}\label{eq:small71}
\|\mathscr{R}\|_{\mathcal{L}(L^2(\partial D), L^2(\partial D))}\leq C e(\tau),
\end{equation}
where
\[
e(\tau)=\begin{cases}
\tau^2\ln\tau\ \ & \mbox{when\ $N=2$},\\
\tau^2\ \ & \mbox{when\ $N=3$}.
\end{cases}
\]
Using (\ref{eq:small70}), the first equation in (\ref{eq:system1}) can be reformulated into
\begin{equation}\label{eq:int1}
\left[\left(\frac{1}{2}I-{K_0}_{[\partial D]}-\mathscr{R}\right)\phi_1(\tau\ \cdot)\right](\tau x')+\left(DL_{[\partial\Omega]}\phi_2(\cdot)\right)(\tau x')=h(\tau x'),\  x'\in\partial D.
\end{equation}
Then using (\ref{eq:small61}), it is directly verified that
\begin{equation}\label{eq:small80}
\left\|\left(DL_{[\partial\Omega]}\phi_2\right)(\tau\ \cdot)\right\|_{L^2(\partial D)}\leq C \left(\tau^{(N-1)/2}\|\phi_1(\cdot)\|_{L^2(\partial D_\tau)}+\tau^N\|\psi\|_{H^{-1/2}(\partial\Omega)}\right).
\end{equation}
Since $I-\frac{1}{2}{K_0}_{[\partial D]}$ is invertible from $L^2(\partial D)$ to $L^2(\partial D)$ (see \cite{Ver}, \cite[\S7.11]{Tay}), by (\ref{eq:small71}), (\ref{eq:int1}), (\ref{eq:small80}) and (\ref{eq:small3}), one can show that
\begin{equation}\label{eq:final1}
\|\phi_1(\tau\ \cdot)\|_{L^2(\partial D)}\leq C\tau\|\psi\|_{H^{-1/2}(\partial\Omega)}.
\end{equation}
Using the relation $\|\phi_1(\cdot)\|_{L^2(\partial D_\tau)}=\tau^{(N-1)/2}\|\phi(\tau\ \cdot)\|_{L^2(\partial D)}$, one further has from (\ref{eq:final1}) that
\begin{equation}\label{eq:final2}
\|\phi_1\|_{L^2(\partial D_\tau)}\leq C\tau^{(N+1)/2}\|\psi\|_{H^{-1/2}(\partial\Omega)}.
\end{equation}
Then, by (\ref{eq:small61}) and (\ref{eq:final2}), we see that
\begin{equation}\label{eq:final3}
\|\phi_2\|_{C(\partial\Omega)}\leq C\tau^N\|\psi\|_{H^{-1/2}(\partial\Omega)}.
\end{equation}
Finally, by the second equation in (\ref{eq:system1}), we have
\begin{equation}\label{eq:final4}
u(x)=2\left\{(DL_{[\partial D_\tau]}\phi_1)(x)-(K_{[\partial\Omega]}\phi_2)(x)+h(x)\right\},\quad x\in\partial \Omega.
\end{equation}
By a similar argument to (\ref{eq:small5n}), one can show that $\|h\|_{C^1(\partial\Omega)}\leq C\tau^N\|\psi\|_{H^{-1/2}(\partial\Omega)}$, which implies
\begin{equation}\label{eq:f1}
\|h\|_{H^{1/2}(\partial\Omega)}\leq C \tau^{N}\|\psi\|_{H^{-1/2}(\partial\Omega)}.
\end{equation}
Since for $x\in\partial\Omega$, $DL_{[\partial D_\tau]}$ has a smooth kernel, by (\ref{eq:final2}) and Schwartz inequality, it is straightforward to show that
\begin{equation}\label{eq:f2}
\|DL_{[\partial D_\tau]}\phi_1\|_{H^{1/2}(\partial\Omega)}\leq \widetilde{C} [\mbox{measure}(\partial D_\tau)]^{1/2}\|\phi_1\|_{L^2(\partial D_\tau)}\leq C\tau^N\|\psi\|_{H^{-1/2}(\partial\Omega)}.
\end{equation}
By noting that $K_{[\partial\Omega]}$ is bounded from $L^2(\partial\Omega)$
to $H^1(\partial\Omega)$ (cf. \cite{ColKre}), we see from (\ref{eq:final3}) that
\begin{equation}\label{eq:f3}
\|K_{[\partial\Omega]}\phi_2\|_{H^{1/2}(\partial\Omega)}\leq C \tau^N \|\psi\|_{H^{-1/2}(\partial\Omega)}.
\end{equation}
Combining (\ref{eq:final4})--(\ref{eq:f3}), we have
\[
\|u\|_{H^{1/2}(\partial\Omega)}\leq C\tau^N\|\psi\|_{H^{-1/2}(\partial\Omega)},
\]
which completes the proof.
\end{pf}

\begin{lem}\label{lem:imp1}
Let $B_{r_l}, l=0,1,2$, $D, D_\tau$ be the ones in (\ref{eq:r0}) and (\ref{eq:r12}) and let $W\in H_{loc}^1(\mathbb{R}^N\backslash\overline{D}_\tau)$ be the unique solution
to
\begin{equation}\label{eq:s11}
\begin{cases}
\ \ \Delta W+\omega^2 W=0\quad & \mbox{in\ $\mathbb{R}^N\backslash\overline{D}_\tau$},\\
\ \  \displaystyle{\frac{\partial W}{\partial \nu}=\phi\in H^{-1/2}(\partial D_\tau)}\quad & \mbox{on\ $\partial D_\tau$},\\
\ \  \displaystyle{\lim_{|x|\rightarrow +\infty}|x|^{(N-1)/2}\left\{\frac{\partial W}{\partial |x|}-i\omega W\right\}=0.}
\end{cases}
\end{equation}
There exists a constant $\tau_0>0$ such that for any $\tau<\tau_0$,
\begin{equation}\label{eq:lemimportant1}
\|W\|_{L^2(B_{r_2}\backslash {B}_{r_0})}\leq C \tau^{N-1}\|\phi(\tau\ \cdot)\|_{H^{-3/2}(\partial D)},
\end{equation}
where $C$ is a positive constant dependent only on $D$ and $\tau_0$, $\omega$, $r_0, r_2$, but independent of $\phi$ and $\tau$.
\end{lem}

In order to gain more insights, we first present a proof of Lemma~\ref{lem:imp1} within spherical geometry. That is, we shall first assume that $D_\tau=B_\tau$, the central ball of radius $\tau$ in $\mathbb{R}^N$, and $D=B_1$. We shall make essential use of series representation of the wave field $W$.

\begin{pf}[Proof of Lemma~\ref{lem:imp1}]

Clearly, in order to show (\ref{eq:lemimportant1}), it suffices to prove that there exist two constants $C_{1}, C_{2}$, such that for all $\phi(\tau \cdot) \in H^{-1/2}(\partial B_{1})$ and $\tau < \tau_{0}$,
\begin{equation}
\begin{cases}
\displaystyle{\|W\|_{H^{1/2}(\partial B_{r_{0}})} \leq C_{1} \tau^{N-1} \|\phi(\tau \cdot)\|_{H^{-3/2}(\partial B_{1})}}, \\
\displaystyle{\|W\|_{H^{1/2}(\partial B_{r_{2}})} \leq C_{2} \tau^{N-1} \|\phi(\tau \cdot)\|_{H^{-3/2}(\partial B_{1})} }.
\end{cases}
\end{equation}

We shall make use of eigenvalues and eigenfunctions of the Laplace-Beltrami operator on a sphere $\partial B_{r}$ to define the Sobolev space $H^s(\partial B_r)$ for $r>0$ and $s\in\mathbb{R}$, which we briefly review in the following and we refer to \cite[\S5.4]{Ned} and \cite[\S1.7]{JE} for general discussions. The Laplace-Beltrami operator on a circle $\partial B_{r}$ in two dimensions is
\[
\Delta_{\partial B_{r}} u = \frac{1}{r^2} \frac{\partial^2 u}{\partial \theta^2},
\]
where $(r,\theta)$ is the polar coordinate in $\mathbb{R}^2$ (cf. \cite[page 234]{JRD});
and the Laplace-Beltrami operator on a sphere $\partial B_{r}$ in three dimensions is
\[
\Delta_{\partial B_{r}} u = \frac{1}{r^2 \sin{\theta}}\left[\frac{\partial}{\partial \theta}(\sin{\theta}\frac{\partial u}{\partial \theta})+ \frac{1}{\sin{\theta}}\frac{\partial^2 u}{\partial \varphi^2}\right],
\]
where $(r, \theta, \varphi)$ is the spherical coordinates in $\mathbb{R}^3$ (see \cite[Appendix]{TZ} and \cite{JRD}).
By direct calculations, one has
\[
-\Delta_{\partial B_{r}}\frac{e^{in\theta}}{\sqrt{2\pi r}} = \frac{n^2}{r^2} \frac{e^{in\theta}}{ \sqrt{2\pi r}}.
\]
$\{w_{n} := \frac{e^{in\theta}}{\sqrt{2\pi r}}\}_{n= -\infty}^{\infty}$ is an orthonormalized  basis in $L^{2}(\partial B_{r})$ and $\lambda_{n}:={n^2}/{r^2}$ is the eigenvalue corresponding to the eigenfunction $w_{n}$. Suppose $u(x)=\sum_{n= -\infty}^{\infty}c_{n}(r)e^{in\theta} \in H^{s}(\partial B_{r})$, $s \in \mathbb{R}$,
then by \cite{Ned} we know
\begin{equation}\label{eq:norm}
\|u\|_{H^{s}(\partial B_{r})}^2 = \sum_{n = -\infty}^{\infty} (1+ n^2/r^2)^s|c_{n}(r)\sqrt{2\pi r}|^2.
\end{equation}
It is noted that when $r = 1$, this is consistent with the $H^{s}[0, 2\pi]$ presented in \cite[\S 8.1]{Rainer}.
In three dimensions, by straightforward calculations, we have
\[
-\Delta_{\partial B_{r}} \frac{Y_{n}^{m}(\hat{x})}{r} = \frac{1}{r^2}n(n+1)\frac{Y_{n}^{m}(\hat{x})}{r},
\]
where $Y_n^m(\hat{x})$ for $\hat{x}\in\mathbb{S}^2$, $n\in\mathbb{N}\cup\{0\}$ and $m=-n, -(n-1),\ldots, (n-1), n$, are the spherical harmonics.
Hence in 3D, the eigenvalues and eigenfunctions are, respectively, $\lambda_{n}= n(n+1)/r^2$ and $w_{n}^m=\frac{Y_{n}^{m}(\hat{x})}{r}$, $m =-n, -(n-1),..., (n-1),n$. Suppose $u(x)=\sum_{n= 0}^{\infty}\sum_{m = -n}^{n}a_{n}^{m}(r)Y_{n}^{m}(\hat{x})\in H^{s}(\partial B_{r})$, $s \in \mathbb{R}$, then we have
\begin{equation}\label{eq:norm1}
\|u\|_{H^{s}(\partial B_{r})}^2 =\sum_{n= 0}^{\infty}\sum_{m=-n}^{n}\left(1+\frac{n(n+1)}{r^2}\right)^s |a_{n}^{m}(r)r|^2.
\end{equation}

We first consider the two-dimensional case of the lemma.
Suppose $\tau < \tau_{0} < r_{0}/4 $ and $\phi(x) = \phi(\tau x') = \sum_{n = -\infty}^{\infty}c_{n}(\tau) e^{in \theta}$, $x  = \tau x', x' \in \partial B_{1}$. Since $\phi(\tau x')\in H^{-1/2}(\partial B_{1})$, one first sees from (\ref{eq:norm})
\[
\|\phi(\tau \cdot)\|_{H^{-3/2}(\partial B_{1})} = \left\{\sum_{n = -\infty}^{\infty}(1+n^2)^{-3/2} 2\pi|c_{n}|^2\right\}^{1/2} < +\infty.
\]
Suppose that the solution is of the form $W(x)=\sum_{n=-\infty}^{\infty}a_{n}H_{n}^{(1)}(\omega |x|)e^{in \theta}$, $x \in\mathbb{R}^2\backslash\overline{B}_\tau$.  According to the PDE system (\ref{eq:s11}), we have $a_{n} = \frac{c_{n}}{\omega{H_{n}^{(1)}}'(\omega \tau)}$ and
\[
W(x)|_{\partial B_{r_{0}}}=  \sum_{n = -\infty}^{\infty} \frac{c_{n}}{\omega{H_{n}^{(1)}}'(\omega\tau)} H_{n}^{(1)}(\omega r_{0}) \sqrt{2 \pi r_{0}} \frac{e^{in\theta}}{\sqrt{2 \pi r_{0}}}.
\]
Hence, by direct calculations
\begin{equation}\label{eq:norm2}
\begin{split}
&\|W(x)\|_{H^{1/2}(\partial B_{r_{0}})} ^2 = \sum_{n=-\infty}^{\infty} (1+n^2/r_{0}^2)^{1/2}\left|\frac{c_{n}H_{n}^{(1)}(\omega r_{0})}{\omega{H_{n}^{(1)}}'(\omega\tau)} \sqrt{2 \pi r_{0}}\right|^{2}\\
\leq& \left\{\sum_{n = -\infty}^{\infty}(1+n^2)^{-3/2} 2\pi|c_{n}|^2\right\}\left\{(1+1/r_{0}^2)^{1/2}\sum_{n = -\infty}^{\infty}(1+n^2)^{2}\left|\frac{ H_{n}^{(1)}(\omega r_{0})}{\omega{H_{n}^{(1)}}'(\omega \tau)} \sqrt{r_{0}}\right|^{2}\right\},
\end{split}
\end{equation}
where we made use of the fact that $(1+ n^2/r_{0}^2) \leq (1+1/r_{0}^2)(1+ n^2)$.
By (\ref{eq:norm1}) and (\ref{eq:norm2}), it is sufficient for us to study the asymptotic development of
\[
T(\tau):=(1+1/r_{0}^2)^{1/2} \sum_{n = -\infty}^{\infty}(1+n^2)^{2}\left|\frac{H_{n}^{(1)}(\omega r_{0})}{\omega{H_{n}^{(1)}}'(\omega\tau)} \sqrt{r_{0}}\right|^{2}.
\]
Since $H_{-n}^{(1)}(\omega r)=(-1)^nH_{n}^{(1)}(\omega r)$, we only need consider $n\geq 0$ in the above series.
By the asymptotic behaviors of Hankel functions as $\tau \rightarrow +0$ or $n \rightarrow +\infty$ (see \cite{AI} and \cite{LiuIMA}), it can be shown that there exists a positive integer $N_{0}$ such that
\begin{equation}\label{2d:mid}
\begin{cases}
\displaystyle{\frac{H_{0}^{(1)}(\omega r_{0})}{\omega{H_{0}^{(1)}}'(\omega\tau)} \sim -\frac{ i \pi H_{0}^{(1)}(\omega r_{0})}{2}  \tau,  \ n = 0},\\
\displaystyle{\frac{H_{n}^{(1)}(\omega r_{0})}{\omega{H_{n}^{(1)}}'(\omega\tau)} \sim  -\frac{i \pi H_{n}^{(1)}(\omega r_{0})}{2^n n! \omega} (\omega\tau)^{n+1},  \ \ n < N_{0},} \\
\displaystyle{\frac{H_{n}^{(1)}(\omega r_{0})}{\omega{H_{n}^{(1)}}'(\omega\tau)} \sim  -\frac{\tau}{n} (\frac{\tau}{r_{0}})^{n},  \ \ n \geq N_{0}.}
\end{cases}
\end{equation}
Furthermore, by our assumption on $\tau$, we have
\begin{equation}\label{2d:mid1}
\sum_{n \geq N_{0}}  \frac{(1+n^2)^{2}}{n^2}\left(\frac{\tau}{r_{0}}\right)^{2n} \leq \sum_{n \geq N_{0}} \frac{(1+n^2)^{2}}{n^2}\frac{1}{4^{2n}} < +\infty.
\end{equation}
Using (\ref{2d:mid}) and (\ref{2d:mid1}), one can show by straightforward calculations that
$T(\tau) \leq C_{1}^{2}\tau^{2}$, where $C_1$ is independent of $\tau$ for $\tau_0$ sufficiently small.
Therefore, we have
\begin{equation}
\|W\|_{H^{1/2}(\partial B_{r_{0}})} \leq C_{1} \tau \|\phi(\tau \cdot)\|_{H^{-3/2}(\partial B_{1})}.
\end{equation}
In like manner, one can show
that there exists $C_{2}$ such that
\begin{equation}
\|W\|_{H^{1/2}(\partial B_{r_{2}})}  \leq C_{2} \tau \|\phi(\tau \cdot)\|_{H^{-3/2}(\partial B_{1})}.
\end{equation}
It is interesting to remark that by similar arguments, one can actually show that for any positive integer $s$, we have
\begin{equation}\label{eq:norm31}
\|W\|_{H^{1/2}(\partial B_{r_{0}})}\leq \widetilde{C}_1 \tau \|\phi(\tau \cdot)\|_{H^{-s/2}(\partial B_{1})},
\end{equation}
and
\begin{equation}\label{eq:norm32}
\|W\|_{H^{1/2}(\partial B_{r_{2}})}  \leq \widetilde{C}_{2} \tau \|\phi(\tau \cdot)\|_{H^{-s/2}(\partial B_{1})}.
\end{equation}
where $\widetilde{C}_1$ and $\widetilde{C}_2$ are positive constants dependent only on $s$, $r_0$, $r_2$, $\tau_0$ and $\omega$, but independent of $\tau$ and $\phi$.

The three dimensional case can be proved similarly.
Suppose $\phi(x)\in H^{-1/2}(\partial B_{\tau})$, and
\[
\phi(x) = \sum_{n = 0}^{\infty}\sum_{m = -n}^{n}c_{n}^{m}(|x|)Y_{n}^{m}.
\]
Then,
\[
\|\phi(\tau \cdot)\|_{H^{-3/2}(\partial B_{1})}^2 = \sum_{n = 0}^{\infty}\sum_{m = -n}^{n}(1+n(n+1))^{-3/2}|c_{n}^{m}(\tau)|^2.
\]
Suppose the solution is of the form
\[
W(x)=\sum_{n = 0}^{\infty}\sum_{m = -n}^{n}a_{n}^{m}h_{n}^{(1)}(\omega|x|)Y_{n}^{m}(\hat{x}).
\]
By using the boundary condition on $\partial B_{\tau}$, we have $a_{n}^{m} = \frac{c_{n}^{m}(\tau)}{\omega{h_{n}^{(1)}}'(\omega \tau)}$. Let $d_{n} = (1+n(n+1))$ in the following calculations. It is noted that $(1+n(n+1)/r_{0}^2) \leq (1+1/r_{0}^{2})d_{n}$. We have
\begin{equation}\label{eq:norm3}
\begin{split}
&\|W\|_{H^{1/2}(\partial B_{r_{0}})}^2 = \sum_{n=0}^{\infty}\sum_{m=-n}^{n}(1+n(n+1)/r_{0}^2)^{1/2}\left|\frac{c_{n}^{m}(\tau)h_{n}^{(1)}(\omega r_{0}) r_{0}}{\omega{h_{n}^{(1)}}'(\omega \tau)}\right|^2\\
\leq & \left\{\sum_{n=0}^{\infty}\sum_{m=-n}^{n}d_{n}^{-3/2}|c_{n}^{m}(\tau)|^2\right\}\left\{(1+1/r_{0}^2)^{1/2}\sum_{n=0}^{\infty}d_{n}^2\left|\frac{h_{n}^{(1)}(\omega r_{0})r_{0}}{\omega{h_{n}^{(1)}}'(\omega \tau)}\right|^2\right\}.
\end{split}
\end{equation}
By the asymptotic properties of the spherical Hankel functions $h_{n}^{(1)}(\omega r)$ and their derivatives as $r \rightarrow +0$ or $n \rightarrow +\infty$ (see \cite{AI} and \cite{LiuIMA}), similar to the two dimensional case, one can show that
\[
(1+1/r_{0}^2)^{1/2}\sum_{n=0}^{\infty}d_{n}^2\left|\frac{h_{n}^{(1)}(\omega r_{0})r_{0}}{\omega{h_{n}^{(1)}}'(\omega\tau)}\right|^2 \leq C_{1}^2 \tau^4.
\]
Therefore, we have
\[
\|W\|_{H^{1/2}(\partial B_{r_{0}})} \leq C_{1} \tau^2 \|\phi(\tau \cdot)\|_{H^{-3/2}(\partial B_{1})}.
\]
In like manner, one can show
\[
\|W\|_{H^{1/2}(\partial B_{r_{2}})} \leq C_{2} \tau^2 \|\phi(\tau \cdot)\|_{H^{-3/2}(\partial B_{1})}.
\]
Similar to the estimates in (\ref{eq:norm31}) and (\ref{eq:norm32}) for the two dimensional case, one can derive more general estimates for the three dimensional case as well.
\end{pf}

Next, we shall present the proof of Lemma~\ref{lem:imp1} within general geometry, which is based on layer-potential techniques.

\begin{pf}[Proof of Lemma~\ref{lem:imp1}]
Let
\[
\varphi=W|_{\partial D_\tau}\in H^{1/2}(\partial D_\tau).
\]
Similar to (\ref{eq:rep2}), we have (cf. \cite{McL})
\begin{equation}\label{eq:green}
W(x)=(DL_{[\partial D_\tau]}\varphi)(x)-(SL_{[\partial D_\tau]}\phi)(x),\quad x\in\mathbb{R}^N\backslash\overline{D}_\tau.
\end{equation}
By the jump properties of layer potential operators (cf. \cite{McL}), we have from (\ref{eq:green}) that
\begin{equation}\label{eq:bi1}
\frac{1}{2}\varphi(x)=(K_{[\partial D_\tau]}\varphi)(x)-(S_{[\partial D_\tau]}\phi)(x),\quad x\in\partial D_\tau.
\end{equation}
Next, we only consider the 3D case, and the 2D case could be shown in a similar manner.

We shall first show that
\begin{equation}\label{eq:cru1}
\|\varphi(\tau\ \cdot)\|_{H^{-1/2}(\partial D)}\leq C\tau\|\phi(\tau\ \cdot)\|_{H^{-3/2}(\partial D)}.
\end{equation}
To that end, we first note that by straightforward calculations
\begin{equation}\label{eq:cru2}
(S_{[\partial D_\tau]}\phi)(\tau x')=\tau({S_0}_{[\partial D]}\phi(\tau\ \cdot))(\tau x')+(\mathscr{S}\phi(\tau\ \cdot))(\tau x'),\quad x'\in\partial D,
\end{equation}
where ${S_0}_{[\partial D]}$ is an integral operator with the kernel given by $G_0(x'-y')$ as the one in (\ref{eq:g0}), and $\mathscr{S}$ is an integral operator with the kernel given by
\begin{equation}\label{eq:residue}
\mathscr{L}(x'-y')=\frac{i\omega}{4\pi}\tau^2+\frac{\tau(i\omega\tau)^2}{2!}\frac{|x'-y'|}{4\pi}+\frac{\tau(i\omega\tau)^3}{3!}\frac{|x'-y'|^2}{4\pi}A(|x'-y'|),
\end{equation}
where $A(t)$ is an (real) analytic function in $t\in\mathbb{R}$. Hence, by the mapping properties presented in \cite[\S 4.3]{Ned}, one has
\begin{equation}\label{eq:cru3}
\|(\mathscr{S}\phi(\tau\ \cdot))(\tau x')\|_{H^{-1/2}(\partial D)}\leq C\tau^2\|\phi(\tau\ \cdot)\|_{H^{-3/2}(\partial D)}.
\end{equation}
Moreover, we know (cf. \cite[\S7.11]{Tay},\cite[\S4.4]{Ned})
\begin{equation}\label{eq:cru4}
\|({S_0}_{[\partial D]}\phi(\tau\ \cdot))(\tau x')\|_{H^{-1/2}(\partial D)}\leq C\|\phi(\tau\ \cdot)\|_{H^{-3/2}(\partial D)}.
\end{equation}

Next, by using the decomposition (\ref{eq:small70}) of $K_{[\partial D_\tau]}$, (\ref{eq:bi1}) can be reformulated as
\begin{equation}\label{eq:cruc}
\begin{split}
&\left[\left(\frac{1}{2}I-{K_0}_{[\partial D]}-\mathscr{R}\right)\varphi(\tau\ \cdot)\right](\tau x')\\
=& \tau({S_0}_{[\partial D]}\phi(\tau\ \cdot))(\tau x')+(\mathscr{S}\phi(\tau\ \cdot))(\tau x') ,\  x'\in\partial D.
\end{split}
\end{equation}
By straightforward asymptotic expansions and also using the mapping properties presented in \cite[\S4.3]{Ned}, one can readily show that
\begin{equation}\label{eq:cru5}
\|\mathscr{R}\|_{\mathcal{L}(H^{-1/2}(\partial D), H^{-1/2}(\partial D))}\leq C\tau^2.
\end{equation}
We shall also make use of the fact that (cf. \cite[\S7.11]{Tay})
\begin{equation}\label{eq:cru6}
\frac{1}{2}I-{K_0}_{[\partial D]}\ \ \mbox{is an isomorphism from $H^{-1/2}(\partial D)$ to $H^{-1/2}(\partial D)$}.
\end{equation}
Hence, by combining (\ref{eq:cru2})--(\ref{eq:cru6}), one readily has (\ref{eq:cru1}).

Finally, by taking $x\in B_{r_2}\backslash B_{r_0}$ in (\ref{eq:green}) and using (\ref{eq:cru1}), we have (\ref{eq:lemimportant1}) by
straightforward verification.

\end{pf}

\begin{lem}\label{lem:small inclusion}
Let $D, D_\tau$ and $\Omega$ be the ones in Lemma~\ref{lem:order n} and let
$P\in H^1(\Omega\backslash\overline{D}_\tau)$ satisfy
\begin{equation}\label{eq:nn1}
\begin{cases}
\ \ \Delta P+\omega^2 P=0\quad & \mbox{in\ $\Omega\backslash\overline{D}_\tau$},\\
\ \ \displaystyle{\frac{\partial P}{\partial \nu}=0}\quad & \mbox{on\ $\partial D_\tau$},\\
\ \ \displaystyle{\frac{\partial P}{\partial \nu}=\varphi\in C(\partial\Omega)}\quad & \mbox{on\ $\partial \Omega$}.
\end{cases}
\end{equation}
Suppose $-\omega^2$ is not an eigenvalue of the Laplacian on $\Omega$ with Neumann boundary condition. Then there exists a constant $\tau_0>0$ such that
for $\tau<\tau_0$, (\ref{eq:nn1}) is uniquely solvable and satisfies
\begin{equation}\label{eq:e1}
\|P\|_{H^{1/2}(\partial\Omega)}\leq C \|\varphi\|_{C(\partial\Omega)},
\end{equation}
where $C$ is a positive constant dependent only on $\tau_0$, $\omega$ and $\Omega$, $D$, but independent of $\tau$ and $\varphi$.
\end{lem}

\begin{pf}
We shall make use of layer potential techniques again to show that lemma. The argument would follow a similar spirit to that for proving Lemma~\ref{lem:order n}, but would be comparatively simpler and we shall only sketch in the following.

Clearly, $P\in C^2(\Omega\backslash\overline{D}_\tau)\cap C(\overline{\Omega}\backslash D_\tau)$ is a strong solution. By letting
\[
p_1=P|_{\partial D_\tau}\quad\mbox{and}\quad p_2=P|_{\partial \Omega},
\]
we have
\begin{equation}\label{eq:rep2n1}
P(x)=(DL_{[\partial D_\tau]}p_1)(x)-(DL_{[\partial\Omega]}p_2)(x)+g(x),
\end{equation}
where
\[
g(x)=\int_{\partial \Omega}G(x-y)\varphi(y)\ d\sigma_y
\]
satisfying
\begin{equation}\label{eq:ff1}
\|g(\tau\ \cdot)\|_{C(\partial D)}\leq C \|\varphi\|_{C(\partial\Omega)}\quad\mbox{and}\quad \|g\|_{C(\partial\Omega)}\leq C\|\varphi\|_{C(\partial\Omega)}.
\end{equation}
By the jump properties of double-layer potential operator, we have from (\ref{eq:rep2n1}) the following system of integral equations for $p_1\in C(\partial D_\tau)$ and $p_2\in C(\partial\Omega)$,
\begin{equation}\label{eq:system1n1}
\begin{cases}
\displaystyle{\frac{1}{2}p_1(x)}=&\displaystyle{(K_{[\partial D_\tau]}p_1)(x)-(DL_{[\partial\Omega]}p_2)(x)+g(x),\quad x\in\partial D_\tau}\\
\displaystyle{\frac{1}{2}p_2(x)}=&\displaystyle{(DL_{[\partial D_\tau]}p_1)(x)-(K_{[\partial\Omega]}p_2)(x)+g(x),\quad x\in\partial \Omega}.
\end{cases}
\end{equation}
By a similar scaling and asymptotic argument to that in the proof of Lemma~\ref{lem:order n}, one can show that
\begin{equation}\label{eq:ff2}
\|p_1\|_{C(\partial D_\tau)}\leq C\|\varphi\|_{C(\partial\Omega)}\quad\mbox{and}\quad \|p_2\|_{C(\partial\Omega)}\leq C\|\varphi\|_{C(\partial\Omega)},
\end{equation}
which in combination with the second equality in (\ref{eq:system1n1}) then implies (\ref{eq:e1}) by direct verification.
\end{pf}

\section{Spherical cloaking device with uniform cloaked contents and sharpness of our estimates}

In this section, we consider our near-cloaking scheme within spherical geometry and uniform cloaked contents. For this special case, we shall assess the cloaking performance, namely Theorem~\ref{thm:main}, and the result illustrates the sharpness of our estimate in Section~\ref{sect:3}.

In the rest of this section, we choose $\Omega$ to be $B_R$, $R>0$, and $\sigma_a'$ to be a scalar constant multiple of the identity matrix, and $q_a'$ to be a positive constant. By a bit abusing of notation, we shall regard $\sigma_a'$ as a scalar constant. In the following, we first consider the two-dimensional case. By transformation acoustics, it is straightforward to show that $\sigma_{a} = \sigma_{a}', q_{a} = q_{a}'/\rho^2$ in $D_{\rho/2}$. Let $\omega_{a} = \omega \sqrt{q_{a}/\sigma_{a}}$ and $\omega_{l} = \omega \sqrt{q_{l}/\sigma_{l}} = \omega \sqrt{1+i} \rho^{-1-\frac{\delta}{2}}$(we choose the branch of $\sqrt{1+i}$  such that $\Im{\sqrt{1+i}}>0$, that is $\sqrt{1+i} = 2^{\frac{1}{4}}e^{i \frac{\pi}{8}}$). Suppose
\[
\psi(x) = \sum_{n=-\infty}^{\infty}\psi_{n}(R)e^{in\theta} \in H^{-1/2}(\partial B_{R}),
\]
and according to our earlier discussion in Section~4,
\begin{equation}\label{normR}
\| \psi\|_{H^{-1/2}(\partial B_{R})}^2 = \sum_{n = -\infty}^{\infty}(1+n^2/R^2)^{-1/2}|\psi_{n}\sqrt{2 \pi R}|^2.
\end{equation}
We assume that the solution of (\ref{eq:virtual wave}) is given by
\begin{equation}\label{virtual-series}
u_\rho(x)=\begin{cases}
\displaystyle{\sum_{n= -\infty}^{\infty} e_{n}J_{n}(\omega_{a}|x|)e^{in \theta}, \ \ x \in {B}_{\rho/2},} \\
\displaystyle{\sum_{n=-\infty}^{\infty}c_{n}J_{n}(\omega_{l}|x|)e^{i n \theta} +  \sum_{n=-\infty}^{\infty}d_{n}H_{n}^{(1)}(\omega_{l}|x|)e^{in \theta},\ \ x \in {B}_{\rho}\backslash \overline{B}_{\rho/2},} \\
\displaystyle{\sum_{n=-\infty}^{\infty}a_{n}J_{n}(\omega|x|)e^{i n \theta} +  \sum_{n=-\infty}^{\infty}b_{n}H_{n}^{(1)}(\omega|x|)e^{in \theta}, \ \ x \in {B}_{R}\backslash \overline{B}_{\rho}.}
\end{cases}
\end{equation}
We shall denote $u_a=u_\rho|_{B_{\rho/2}}$, $u_l=u_\rho|_{B_\rho\backslash\overline{B}_{\rho/2}}$ and $u_R=u|_{B_R\backslash\overline{B}_{\rho}}$.
By the standard transmission conditions on $\partial B_{\rho/2}, \partial B_{\rho}$ and the boundary condition on $\partial B_{R}$, we have
\begin{equation}\label{boundary:virtual}
\begin{cases}
\displaystyle{u_{a}(x) = u_{l}(x), \ \sigma_{a}\frac{\partial u_{a}(x)}{\partial \nu(x)} = \sigma_{l}\frac{\partial u_{l}(x)}{\partial \nu(x)}, \ \ x \in \partial B_{\rho/2},} \\
\displaystyle{u_{l}(x) = u_{R}(x),  \ \sigma_{l}\frac{\partial u_{l}(x)}{\partial \nu(x)} = \frac{\partial u_{R}(x)}{\partial \nu(x)}, \ \ x \in  \partial B_{\rho},} \\
\displaystyle{\frac{\partial u_{R}(x)}{\partial \nu(x)} = \psi(x), \ \ x \in \partial B_{R}.}
\end{cases}
\end{equation}
Plugging the series representations (\ref{virtual-series}) into (\ref{boundary:virtual}), we have the following linear system of equations for the coefficients,
\begin{equation}\label{series:equ}
\begin{cases}
e_{n}J_{n}(\omega_{a} \rho/2) = c_{n}J_{n}(\omega_{l}\rho/2) + d_{n}H_{n}^{(1)}(\omega_{l}\rho/2), \\
\sqrt{\sigma_{a}q_{a}}e_{n}J_{n}'(\omega_{a} \rho/2) = \sqrt{\sigma_{l}q_{l}}[c_{n}J_{n}'(\omega_{l}\rho/2) + d_{n}{H_{n}^{(1)}}'(\omega_{l}\rho/2)], \\
c_{n}J_{n}(\omega_{l}\rho) + d_{n}H_{n}^{(1)}(\omega_{l}\rho) = a_{n}J_{n}(\omega \rho) + b_{n} H_{n}^{(1)}(\omega \rho), \\
 \sqrt{\sigma_{l}q_{l}}[c_{n}J_{n}'(\omega_{l}\rho) + d_{n}{H_{n}^{(1)}}'(\omega_{l}\rho)] = a_{n}J_{n}'(\omega \rho) + b_{n} {H_{n}^{(1)}}'(\omega \rho), \\
 a_{n}\omega J_{n}'(\omega R) + b_{n}\omega {H_{n}^{(1)}}'(\omega R) = \psi_{n}.
\end{cases}
\end{equation}
Letting $A = \sqrt{\frac{\sigma_{a}q_{a}}{\sigma_{l}q_{l}}} = \frac{\sqrt{\sigma_{a}'q_{a}'}}{2^{\frac{1}{4}}e^{i\frac{\pi}{8}}}\rho^{-2-\frac{\delta}{2}}$, from the first two equations of (\ref{series:equ}) we have
\begin{equation}\label{eq:dncn}
\begin{cases}
d_{n} = -\frac{J_{n}(\omega_{l}\rho/2)}{H_{n}^{(1)}(\omega_{l}\rho/2)}c_{n}\ \ \mbox{if} \ \ J_{n}(\omega_{a}\rho/2) = 0, \\
d_{n} = -\frac{J_{n}'(\omega_{l}\rho/2)-AJ_{n}(\omega_{l}\rho/2)\frac{J_{n}'(\omega_{a}\rho/2)}{J_{n}(\omega_{a}\rho/2)}}{{H_{n}^{(1)}}'(\omega_{l}\rho/2)-A H_{n}^{(1)}(\omega_{l}\rho/2) \frac{J_{n}'(\omega_{a}\rho/2)}{J_{n}(\omega_{a}\rho/2)}}c_{n}\ \ \mbox{if} \ \ J_{n}(\omega_{a}\rho/2) \neq 0.
\end{cases}
\end{equation}
Denoting the expressions before $c_n$ in (\ref{eq:dncn}) by $\Upsilon_n$, namely $d_{n}: =\Upsilon_{n}c_{n}$, and substituting $d_n$ into the third and fourth equations of (\ref{series:equ}), we have by straightforward calculations
\begin{equation}\label{bn:an}
b_{n} = -\frac{2^{\frac{1}{4}}e^{i\frac{\pi}{8}}  \rho^{1+ \frac{\delta}{2}} \frac{J_{n}'(\omega_{l}\rho) + \Upsilon_{n}{H_{n}^{(1)}}'(\omega_{l}\rho)}{J_{n}(\omega_{l}\rho) + \Upsilon_{n} H_{n}^{(1)}(\omega_{l}\rho)}J_{n}(\omega \rho) - J_{n}'(\omega\rho)}{2^{\frac{1}{4}}e^{i \frac{\pi}{8}}\rho^{1+ \frac{\delta}{2}} \frac{J_{n}'(\omega_{l}\rho) + \Upsilon_{n}{H_{n}^{(1)}}'(\omega_{l}\rho)}{J_{n}(\omega_{l}\rho) + \Upsilon_{n} H_{n}^{(1)}(\omega_{l}\rho)}H_{n}^{(1)}(\omega\rho) -{H_{n}^{(1)}}'(\omega\rho)}a_{n}.
\end{equation}
Let $\Gamma_n$ denote the expression before $a_n$ in (\ref{bn:an}), namely $b_{n}: = \Gamma_{n} a_{n}$, and
\begin{equation}\label{eq:H}
\mathcal{H}_{n}(\rho) = \frac{J_{n}'(\omega_{l}\rho) + \Upsilon_{n}{H_{n}^{(1)}}'(\omega_{l}\rho)}{J_{n}(\omega_{l}\rho) + \Upsilon_{n} H_{n}^{(1)}(\omega_{l}\rho)}.
\end{equation}
Plugging (\ref{bn:an}) into the last equation in (\ref{virtual-series}), we have
\begin{equation}\label{eq:ur}
u_{R}(x) = \sum_{n = -\infty}^{\infty} \frac{\psi_{n}[J_{n}(\omega R) + \Gamma_{n}H_{n}^{(1)}(\omega R)]}{\omega[J_{n}'(\omega R) + \Gamma_{n} {H_{n}^{(1)}}'(\omega R)]}e^{i n \theta}, \ \ \ x \in \partial B_{R},
\end{equation}
whereas the ``free space" solution $u_{0}(x) \in H^{1}(\Omega)$ of (\ref{eq:free space}) is
\begin{equation}
u_{0}(x) = \sum_{n = -\infty}^{\infty} \frac{\psi_{n} J_{n}(\omega |x|)}{\omega J_{n}'(\omega R)}e^{i n\theta}, \ \  \ x \in B_{R}.
\end{equation}
Hence,
\begin{equation}
\left[u_{\rho}(x)-u_{0}(x)\right]|_{\partial B_{R}}  = \sum_{n = -\infty}^{\infty} \frac{\psi_{n} J_{n}(\omega R)}{\omega J_{n}'(\omega R)}\left[ \frac{\Gamma_{n}\left[\frac{H_{n}^{(1)}(\omega R)}{J_{n}(\omega R)} -\frac{{H_{n}^{(1)}}'(\omega R)}{J_{n}'(\omega R)}\right]}{1+ \frac{\Gamma_{n}{H_{n}^{(1)}}'(\omega R)}{J_{n}'(\omega R)}}\right]e^{i n \theta},
\end{equation}
and therefore
\begin{equation}\label{u:rho}
\begin{split}
&\|u_{\rho} -u_{0}\|_{H^{1/2}(\partial B_{R})}^{2} = \sum_{n = -\infty}^{\infty} \left(1+ \frac{n^{2}}{R^{2}}\right)^{1/2} \left|\frac{\psi_{n} J_{n}(\omega R)}{\omega J_{n}'(\omega R)}\widetilde{h}_{n}\sqrt{2\pi R}\right|^2\\
\leq & \left\{\sum_{n = -\infty}^{\infty} \left(1+ \frac{n^{2}}{R^{2}}\right)^{-1/2}|\psi_{n}\sqrt{2\pi R}|^{2} \right\} \left\{\sum_{n = -\infty}^{\infty} \left(1+ \frac{n^{2}}{R^{2}}\right) \left|\frac{J_{n}(\omega R)}{\omega J_{n}'(\omega R)}\widetilde{h}_{n}\right|^2\right\}\\
\leq & \left\{\sum_{n = -\infty}^{\infty} \left(1+ \frac{n^{2}}{R^{2}}\right) \left|\frac{J_{n}(\omega R)}{\omega J_{n}'(\omega R)}\widetilde{h}_{n}\right|^2\right\}\|\psi\|_{H^{-1/2}(\partial B_R)}^2
\end{split}
\end{equation}
where
\[
\widetilde{h}_{n} = \left|\frac{\Gamma_{n}\left[\frac{H_{n}^{(1)}(\omega R)}{J_{n}(\omega R)} -\frac{{H_{n}^{(1)}}'(\omega R)}{J_{n}'(\omega R)}\right]}{1+ \frac{\Delta_{n}{H_{n}^{(1)}}'(\omega R)}{J_{n}'(\omega R)}}\right|.
\]

Since $H_{-n}^{(1)}(\omega r)=(-1)^nH_{n}^{(1)}(\omega r)$, we only need consider $n\geq 0$ in estimating the series in the last inequality in (\ref{u:rho}).
Using the asymptotic behaviors of $J_{n}(z)$, $H_{n}^{(1)}(z)$, $J_{n}'(z)$, ${H_{n}^{(1)}}'(z)$ as both $\Im{z}$ and $\Re{z}$ tend to $+\infty$ (cf. \cite {AI}, \cite{LiLiuSun}), one can show
\begin{equation}\label{asym:H}
\mathcal{H}_{n}(\rho) \sim \frac{J_{n}'(\omega_{l}\rho)}{J_{n}(\omega_{l}\rho)} \sim -e^{i\pi/2} +¡¡\mathcal{O}(n\rho^{\frac{\delta}{2}}).
\end{equation}
which together with the asymptotic behaviors of $J_{n}(\omega \rho)$, $H_{n}^{(1)}(\omega \rho)$, $J_{n}'(\omega \rho)$, ${H_{n}^{(1)}}'(\omega \rho)$ as $\rho \rightarrow +0$ (cf. \cite{AI,LiuIMA}), one can further show
\begin{equation}\label{deltan:jn}
\begin{cases}
\Gamma_{0} \sim \frac{\rho^2 \pi \omega i}{2}[2^{\frac{1}{4}} e^{i \frac{5\pi}{8}}\rho^{\frac{\delta}{2}} - \frac{\omega}{2}] , \ \ n = 0,\\
\Gamma_{n} \sim  \pi i(\omega \rho)^{2n} [2^{\frac{5}{4}} e^{i \frac{5\pi}{8}}\omega \rho^{2+\frac{\delta}{2}}+n]/(2^n n!)^2 , \ \  n \geq 1.
\end{cases}
\end{equation}
Then using the estimates in (\ref{deltan:jn}), together with the use of the asymptotic developments of the Bessel and Hankel functions for large $n$ (cf. \cite{AI}), one can verify that there exists a sufficiently large integer $N_{1}$ such that
\begin{equation}\label{deltan:jn1}
\begin{cases}
 \widetilde{h}_{0} \sim \frac{\rho^2 \pi \omega i}{2}[2^{\frac{1}{4}} e^{i \frac{5\pi}{8}}\rho^{\frac{\delta}{2}} - \frac{\omega}{2}][\frac{H_{0}^{(1)}(\omega R)}{J_{0}(\omega R)} - \frac{{H_{0}^{(1)}}'(\omega R)}{J_{0}'(\omega R)}], \ \ n = 0,\\
 \widetilde{h}_{n} \sim  \pi i(\omega \rho)^{2n} \frac{n}{(2^n n!)^2} [\frac{H_{n}^{(1)}(\omega R) }{J_{n}(\omega R)} - \frac{{H_{n}^{(1)}}'(\omega R)}{J_{n}'(\omega R)}], \ \ 1 \leq n \leq N_{1}, \\
 \widetilde{h}_{n} \sim  \frac{2}{n}(\frac{\rho}{R})^{2n}[2^{\frac{5}{4}} e^{i \frac{5\pi}{8}}\omega \rho^{2+\frac{\delta}{2}}+n], \ \ n > N_{1}.
\end{cases}
\end{equation}
Hence from (\ref{deltan:jn1}), we readily see that there exists a constant $C_{1}$ independent of $\rho$ for $\rho$ sufficiently small such that
\begin{equation}\label{optimal}
|\widetilde{h}_{n}| \leq C_{1}\rho^{2}, \ \ n \leq N_{1},
\end{equation}\label{n:bigbig}
and for $n > N_{1}$
\begin{equation}
|\widetilde{h}_{n}| \leq  \frac{8}{n}\left(\frac{\rho}{R}\right)^{2n}[2^{\frac{5}{4}} \omega \rho^{2+\frac{\delta}{2}}+n].
\end{equation}
Here it is emphasized that due to the asymptotic developments of $\widetilde{h}_0$ and $\widetilde{h}_1$, (\ref{optimal}) is the best estimate one could achieve, namely $C_1\rho^2$ could not be improved.
Now, using(\ref{optimal}), we see that
\begin{equation}\label{optimal1}
\sum_{n =0}^{N_{1}}\left(1+ \frac{n^{2}}{R^{2}}\right) \left|\frac{J_{n}(\omega R)}{\omega J_{n}'(\omega R)}\widetilde{h}_{n}\right|^2 \leq C_{2} \rho^{4}.
\end{equation}
Let $N_1$ be sufficiently large such that $|\frac{J_{n}(\omega R)}{\omega J_{n}'(\omega R)}| < 1$ for $n>N_1$, then for $\rho<\min\{R/4, 1\}$
\begin{equation}\label{less:infty}
\begin{split}
&\sum_{n > N_{1}}\left(1+ \frac{n^{2}}{R^{2}}\right) \left|\frac{J_{n}(\omega R)}{\omega J_{n}'(\omega R)}\widetilde{h}_{n}\right|^2\\
\leq& \frac{\rho^{4}}{R^{4}}\sum_{n > N_{1}}\left(1+ \frac{n^{2}}{R^{2}}\right) \left|\frac{8}{n}\left(\frac{\rho}{R}\right)^{2(n-1)}\left[2^{\frac{5}{4}} \omega \rho^{2+\frac{\delta}{2}}+n\right]\right|^{2}<C_3\rho^4.
\end{split}
\end{equation}
Combining (\ref{u:rho}), (\ref{optimal1}) and (\ref{less:infty}),  we have
\begin{equation}\label{eq:sharp2d}
\|u_{\rho} -u_{0}\|_{H^{1/2}(\partial B_{R})} \leq C \rho^{2} \|\psi\|_{H^{-1/2}(\partial B_{R})}.
\end{equation}
Moreover, from the optimality of the estimate (\ref{optimal}), we readily see the sharpness of (\ref{eq:sharp2d}).

Next, we shall investigate the asymptotic behavior of the boundary condition on $\partial B_\rho^+$, namely, $\frac{\partial u_R^+}{\partial\nu}|_{\partial B_\rho}$. Since
\[
\frac{\partial u_{R}^{+}}{\partial \nu}\bigg|_{\partial B_{\rho}} = \sum_{n=-\infty}^{\infty}\omega l_{n}e^{i n\theta},
\]
where
\[
l_{n}:=a_{n}J_{n}'(\omega \rho) + b_{n}{H_{n}^{(1)}}'(\omega \rho)
\]
and from the last equation of (\ref{series:equ})
\begin{equation}
a_{n} = \frac{\psi_{n}}{\omega[J_{n}'(\omega R)+ \Gamma_{n}{H_{n}^{(1)}}'(\omega R)]}.
\end{equation}
Hence, we only need study the asymptotic behavior of $l_n$. By direct calculations, we have
\[
l_{n}= \frac{2^{\frac{1}{4}}e^{i\frac{\pi}{8}}\rho^{1+\frac{\delta}{2}}\mathcal{H}_{n}(\rho)\psi_{n}}{\omega  [ J_{n}'(\omega R)+\Gamma_{n} {H_{n}^{(1)}}'(\omega R)] } \frac{J_{n}'(\omega \rho)H_{n}^{(1)}(\omega \rho)-J_{n}(\omega \rho){H_{n}^{(1)}}'(\omega \rho)}{[2^{\frac{1}{4}}e^{i\frac{\pi}{8}}\rho^{1+\frac{\delta}{2}}\mathcal{H}_{n}(\rho)H_{n}^{(1)}(\omega \rho) - {H_{n}^{(1)}}'(\omega \rho)]}.
\]
Using the asymptotic behavior of $\mathcal{H}_{n}(\rho)$ in (\ref{asym:H}), $\Gamma_{n}$ in (\ref{deltan:jn}), and the Wronskian $J_{n}(t)Y_{n}'(t) - J_{n}'(t)Y_{n}(t) = \frac{2}{\pi t}$ (cf. \cite{ColKre}), one can show that there exists a sufficiently large integer $N_{2}$ such that
\begin{equation}\label{eq:aa1}
\begin{cases}
 l_{0}\sim -2^{\frac{1}{4}}e^{i\frac{5\pi}{8}} \rho^{1+\frac{\delta}{2}} \frac{\psi_{0}}{\omega J_{0}'(\omega R)}, \ \ n =0,\\
 l_{n}\sim -\frac{\psi_{n}2^{\frac{1}{4}}e^{i\frac{5\pi}{8}} \rho^{1+\frac{\delta}{2}}}{\omega J_{n}'(\omega R)} \frac{2(\omega \rho)^{n}}{2^{n}n!}= \mathcal{O}(\rho^{n+1+\delta/2}), \ \ 1 \leq n \leq N_{2},\\
l_{n} \sim  -\frac{2\psi_{n}2^{\frac{1}{4}}e^{i \frac{5\pi}{8}} \rho^{1+\frac{\delta}{2}}R}{n}
(\frac{\rho}{R})^{n}, \ \ n > N_{2}.
\end{cases}
\end{equation}
Using (\ref{eq:aa1}) and a similar argument to that for the proof of Lemma \ref{lem:imp1}, we can show that
\begin{equation}\label{eq:sh1}
\left\|\frac{\partial u_{R}^{+}}{\partial \nu}(\rho\ \cdot)\right\|_{H^{-1/2}(\partial B_{1})} \leq C \rho^{1+\frac{\delta}{2}}\|\psi\|_{H^{-1/2}(\partial B_{R})},
\end{equation}
where $C$ is independent of $\rho$, $\delta$ and $\psi$. Hence, we readily see that as $\delta\rightarrow +\infty$, the lossy layer $\{B_\rho\backslash\overline{B}_{\rho/2}; \sigma_l, q_l\}$ converges to a {\it sound-hard} layer; that is, the normal velocity of the wave filed would vanish on the exterior of the layer, namely
\begin{equation}\label{eq:sh2}
\left\|\frac{\partial u_{R}^{+}}{\partial \nu}(\rho\ \cdot)\right\|_{H^{-1/2}(\partial B_{1})}\rightarrow 0\quad \mbox{as\ $\delta\rightarrow+\infty$}.
\end{equation}
On the other hand, the sound-hard layer lining is considered in \cite{LiLiuSun}, and it is shown that one could achieve optimally $\rho^2$ within the ideal cloaking for the regularized cloaking construction. Hence, (\ref{eq:sh1}) and (\ref{eq:sh2}) also partly illustrate the sharpness of our estimates.

The three-dimensional case could be treated similarly, which we only sketch in the following. Let
\[
\psi(x)|_{\partial B_{R}} = \sum_{n = 0}^{\infty} \sum_{m=-n}^{n}\psi_{n}^{m}Y_{n}^{m}(\hat{x}) \in H^{-1/2}(\partial B_{R}),
\]
with
\begin{equation}
\|\psi(x)\|_{H^{-1/2}(\partial B_{R})}^2 =  \sum_{n = 0}^{\infty} \sum_{m=-n}^{n}(1+n(n+1)/R^2)^{-1/2} |\psi_{n}^{m}R|^2 < +\infty.
\end{equation}
Noting $\sigma_{a} = \sigma_{a}'/\rho$, $q_{a} = q_{a}'/\rho^2$
in $B_{\rho/2}$, similar to (\ref{virtual-series}) for the 2D case, the wave fields in the separated domains could be represented as follows
\begin{equation}\label{eq:3d}
\begin{split}
u_{a}(x) =&\sum_{n = 0}^{\infty} \sum_{m = -n}^{n} e_{n}^{m}j_{n}(\omega_{a}|x|)Y_{n}^{m}(\hat{x}), \\
u_{l}(x) =& \sum_{n = 0}^{\infty} \sum_{m = -n}^{n} c_{n}^{m}j_{n}(\omega_{l}|x|)Y_{n}^{m}(\hat{x}) + \sum_{n = 0}^{\infty} \sum_{m = -n}^{n} d_{n}^{m}h_{n}^{(1)}(\omega_{l}|x|)Y_{n}^{m}(\hat{x}), \\
u_{R}(x) =& \sum_{n = 0}^{\infty} \sum_{m = -n}^{n} a_{n}^{m} j_{n}(\omega |x|)Y_{n}^{m}(\hat{x}) + \sum_{n = 0}^{\infty} \sum_{m = -n}^{n} b_{n}^{m}h_{n}^{(1)}(\omega|x|)Y_{n}^{m}(\hat{x}).
\end{split}
\end{equation}
Similar to (\ref{boundary:virtual}), using the standard transmission conditions and the boundary condition, one could derive the following linear system of equations for the coefficients
\begin{equation}\label{3d:series1}
\begin{cases}
e_{n}^{m}j_{n}(\omega_{a}\rho/2) = c_{n}^{m}j_{n}(\omega_{l} \rho/2) + d_{n}^{m}h_{n}^{(1)}(\omega_{l} \rho/2) \\
\sqrt{\sigma_{a}q_{a}}e_{n}^{m}j_{n}'(\omega_{a}\rho/2) = \sqrt{\sigma_{l}q_{l}}[c_{n}^{n}j_{n}'(\omega_{l} \rho/2) + d_{n}^{m}{h_{n}^{(1)}}'(\omega_{l} \rho/2)]\\
c_{n}^{m}j_{n}(\omega_{l} \rho) + d_{n}^{m}h_{n}^{(1)}(\omega_{l} \rho) = a_{n}^{m}j_{n}(\omega \rho) + b_{n}^{m}h_{n}^{(1)}(\omega \rho) \\
\sqrt{\sigma_{l}q_{l}}[c_{n}^{m}j_{n}'(\omega_{l} \rho) + d_{n}^{m}{h_{n}^{(1)}}'(\omega_{l} \rho)] = a_{n}^{m}j_{n}'(\omega \rho) + b_{n}^{m}{h_{n}^{(1)}}'(\omega \rho), \\
\omega a_{n}^{m}j_{n}'(\omega R) + \omega b_{n}^{m}{h_{n}^{(1)}}'(\omega R) = \psi_{n}^{m}.
\end{cases}
\end{equation}
Letting $\widetilde{A} = \sqrt{\frac{\sigma_{a}q_{a}}{\sigma_{l}q_{l}}} = \frac{\sqrt{\sigma_{a}'q_{a}'}}{2^{\frac{1}{4}}e^{i\frac{\pi}{8}}}\rho^{-\frac{5}{2}-\frac{\delta}{2}}$ and solving (\ref{3d:series1}), one has
\begin{equation}\label{bnm:anm}
b_{n}^{m} = -\frac{2^{\frac{1}{4}}e^{i\frac{\pi}{8}} \rho^{1+ \frac{\delta}{2}} \frac{j_{n}'(\omega_{l}\rho) + \widetilde{\Upsilon}_{n}{h_{n}^{(1)}}'(\omega_{l}\rho)}{j_{n}(\omega_{l}\rho) + \widetilde{\Upsilon}_{n} h_{n}^{(1)}(\omega_{l}\rho)}j_{n}(\omega \rho) - j_{n}'(\omega\rho)}{2^{\frac{1}{4}}e^{i\frac{\pi}{8}} \rho^{1+ \frac{\delta}{2}}\frac{j_{n}'(\omega_{l}\rho) + \widetilde{\Upsilon}_{n}{h_{n}^{(1)}}'(\omega_{l}\rho)}{j_{n}(\omega_{l}\rho) + \widetilde{\Upsilon}_{n} h_{n}^{(1)}(\omega_{l}\rho)}h_{n}^{(1)}(\omega\rho) -{h_{n}^{(1)}}'(\omega\rho)}a_{n}^{m},
\end{equation}
where
\begin{equation}
\widetilde{\Upsilon}_{n}:=\begin{cases}
-\frac{j_{n}(\omega_{l}\rho/2)}{h_{n}^{(1)}(\omega_{l}\rho/2)} \ \ \mbox{if}\ \ j_{n}(\omega_{a}\rho/2) = 0, \\
-\frac{j_{n}'(\omega_{l}\rho/2)-\widetilde{A}j_{n}(\omega_{l}\rho/2)\frac{j_{n}'(\omega_{a}\rho/2)}{j_{n}(\omega_{a}\rho/2)}}{{h_{n}^{(1)}}'(\omega_{l}\rho/2)-\widetilde{A} h_{n}^{(1)}(\omega_{l}\rho/2) \frac{j_{n}'(\omega_{a}\rho/2)}{j_{n}(\omega_{a}\rho/2)}} \ \ \mbox{if}\ \ j_{n}(\omega_{a}\rho/2) \neq 0.
\end{cases}
\end{equation}
Let $\widetilde{\Gamma}_n$ denote the expression before $a_n^m$ in (\ref{bnm:anm}). Then
\begin{equation}
[u_{\rho}(x)-u_{0}(x)]|_{\partial B_{R}}  = \sum_{n = 0}^{\infty}\sum_{m=-n}^{n} \frac{\psi_{n}^{m} j_{n}(\omega R)}{\omega j_{n}'(\omega R)}\left\{ \frac{\widetilde{\Gamma}_{n}\left[\frac{h_{n}^{(1)}(\omega R)}{j_{n}(\omega R)} -\frac{{h_{n}^{(1)}}'(\omega R)}{j_{n}'(\omega R)}\right]}{1+ \frac{\widetilde{\Gamma}_{n}{h_{n}^{(1)}}'(\omega R)}{j_{n}'(\omega R)}}\right\}Y_{n}^{m}(\hat{x}).
\end{equation}
Let
\[
\widetilde{g}_{n} =  \widetilde{\Gamma}_{n}\left[\frac{h_{n}^{(1)}(\omega R)}{j_{n}(\omega R)} -\frac{{h_{n}^{(1)}}'(\omega R)}{j_{n}'(\omega R)}\right]\bigg/ \left[1+ \frac{\widetilde{\Gamma}_{n}{h_{n}^{(1)}}'(\omega R)}{j_{n}'(\omega R)}\right]
.\]
Then
\begin{equation}\label{eq:3dd}
\begin{split}
&\|u_{\rho}(x)-u_{0}(x)\|_{H^{1/2}(\partial B_{R})}^2 = \sum_{n = 0}^{\infty}\sum_{m=-n}^{n}  \sqrt{1+\frac{n(n+1)}{R^2}} \left| \frac{\psi_{n}^{m} j_{n}(\omega R)}{\omega j_{n}'(\omega R)}\widetilde{g}_{n}R\right|^2\\
\leq &\left\{\sum_{n = 0}^{\infty}\sum_{m=-n}^{n}\frac{1}{\sqrt{1+\frac{n(n+1)}{R^2}}} \left|\psi_{n}^{m} R\right|^2\right\}\left\{\sum_{n=0}^{\infty}(1+\frac{n(n+1)}{R^2})\left|\frac{j_{n}(\omega R) \widetilde{g}_{n}}{\omega j_{n}'(\omega R)}\right|^2\right\}.
\end{split}
\end{equation}
By similar asymptotic analyses to the 2D case, one can show that there exists a sufficiently large integer $N_3$ such that
\begin{equation}\label{asym:3d}
\begin{cases}
\widetilde{g}_{0} \sim i\rho^{3}[2^{\frac{1}{4}}e^{i\frac{5\pi}{8}}\rho^{\frac{\delta}{2}}\omega^2 - \frac{\omega^3}{3}][\frac{h_{0}^{(1)}(\omega R)}{j_{0}(\omega R)} -\frac{{h_{0}^{(1)}}'(\omega R)}{j_{0}'(\omega R)}], \ n=0, \\
\widetilde{g}_{n} \sim i\frac{(\omega \rho)^{2n+1}n(n!2^n)^2}{(n+1)(2n)!(2n+1)!}[\frac{h_{n}^{(1)}(\omega R)}{j_{n}(\omega R)} -\frac{{h_{n}^{(1)}}'(\omega R)}{j_{n}'(\omega R)}] , \ 1 \leq n < N_{3},\\
\widetilde{g}_{n} \sim (\frac{\rho}{R})^{2n+1}\frac{2n+1}{n(n+1)}[n+ \frac{2n+1}{n+1}2^{\frac{1}{4}}e^{i\frac{5\pi}{8}} \omega \rho^{2+ \frac{\delta}{2}}], \ n \geq N_{3}.
\end{cases}
\end{equation}
By a similar argument to the 2D case, applying (\ref{asym:3d}) to the estimation of (\ref{eq:3dd}), one can show that
\begin{equation}\label{eq:sh3}
\|u_{\rho}(x)-u_{0}(x)\|_{H^{1/2}(\partial B_{R})} \leq C \rho^3\|\psi\|_{H^{-1/2}(\partial B_{R})},
\end{equation}
and moreover, the estimate is optimal. Furthermore, one could also show in this 3D case
\begin{equation}\label{eq:sh3}
\left\|\frac{\partial u_{R}^{+}}{\partial \nu}(\rho \cdot )\right\|_{H^{-1/2}(\partial B_{1})} \leq C \rho^{1+\frac{\delta}{2}} \|\psi\|_{H^{-1/2}(\partial B_{R})}.
\end{equation}
That is, we also have that the lossy layer would converge to a sound-hard layer in the limiting case as $\delta\rightarrow+\infty$.

\section{Discussion}

In this work, we consider a novel near-cloaking scheme by employing a well-designed lossy layer between the cloaked region and the cloaking region. The study follows the spirit of the one developed in \cite{KOVW}. However, in \cite{KOVW} the authors rely on a lossy layer with a large lossy parameter for the successful near-cloaking construction, whereas we rely on a lossy layer with a large density parameter. They are of different physical and mathematical nature.
As was discussed earlier in Introduction, the lossy layer proposed in \cite{KOVW} is a finite realization of a sound-soft layer, whereas the FSH layer in the current work is a finite realization of a sound-hard layer. This is confirmed by (\ref{eq:sh1}) and (\ref{eq:sh3}) derived in Section 5, which indicates that as $\delta\rightarrow+\infty$ the FSH layer converges to a sound-hard layer. Moreover, we have the following result which further supports our such observation.
\begin{thm}\label{thm:sh}
Suppose $-\omega^2$ is not an eigenvalue of the Laplacian on $\Omega\backslash\overline{D}$ with Neumann boundary condition. Let $u_{sh}\in H^1(\Omega\backslash\overline{D})$ be the unique solution of
\begin{equation}\label{eq:soundhard}
\begin{cases}
& \displaystyle{\sum_{i,j=1}^N\frac{\partial}{\partial x_i}\left((\sigma_c^\rho)^{ij}(x)\frac{\partial u_{sh}}{\partial x_j}\right)+\omega^2 q_c^\rho(x) u_{sh}=0\quad\mbox{in\ \ $\Omega\backslash\overline{D}$},}\\
& \displaystyle{\sum_{i,j=1}^N\nu_i(\sigma_c^\rho)^{ij}\frac{\partial u_{sh}}{\partial x_j}=\psi\in H^{-1/2}(\partial\Omega)\quad\mbox{on\ \ $\partial\Omega$},}\\
&\displaystyle{\sum_{i,j=1}^N\tilde{\nu}_i(\sigma_c^\rho)^{ij}\frac{\partial u_{sh}}{\partial x_j}=0\quad\mbox{on\ \ $\partial D$},}
\end{cases}
\end{equation}
where $\{\Omega; \sigma_c^\rho, q_c^\rho\}$ is the medium in (\ref{eq:cloaking medium}) and, $\nu=(\nu_i)_{i=1}^N$ and $\tilde{\nu}=(\tilde{\nu}_i)_{i=1}^N$
are the exterior unit normals to $\partial D$ and $\partial\Omega$, respectively. That is, $u_{sh}$ is the solution corresponding to a sound-hard obstacle $D$ buried in the medium $\{\Omega; \sigma_c^\rho, q_c^\rho\}$. Let $u\in H^1(\Omega)$ be the solution corresponding to the cloaking device, namely,
\begin{equation}\label{eq:fsh}
\begin{cases}
& \displaystyle{\sum_{i,j=1}^N\frac{\partial}{\partial x_i}\left(\sigma^{ij}(x)\frac{\partial u}{\partial x_j}\right)+\omega^2 q(x) u=0\quad\mbox{in\ \ $\Omega$},}\\
& \displaystyle{\sum_{i,j=1}^N\nu_i\sigma^{ij}\frac{\partial u}{\partial x_j}=\psi\in H^{-1/2}(\partial\Omega)\quad\mbox{on\ \ $\partial\Omega$}.}
\end{cases}
\end{equation}
Then for sufficiently small $\rho>0$, we have
\begin{equation}\label{eq:c1}
\|u_{sh}-u\|_{H^{1/2}(\partial\Omega)}\leq C\rho^N\|\psi\|_{H^{-1/2}(\partial\Omega)},
\end{equation}
where $C$ is a constant independent of $\rho$ and $\psi$.
\end{thm}

\begin{pf}
Let
\[
\tilde{u}_{sh}=F^*u_{sh}\quad\mbox{and}\quad \tilde{u}=F^* u,
\]
then by transformation acoustics, it is readily seen that $\tilde{u}$ is exactly $u_\rho$ in (\ref{eq:virtual wave}) and $\tilde{u}_{sh}$
satisfies
\begin{equation}\label{eq:soundhard1}
\begin{cases}
&\displaystyle{ \Delta \tilde{u}_{sh}+\omega^2 \tilde{u}_{sh}=0\quad\mbox{in\ \ $\Omega\backslash\overline{D}_\rho$},}\\
&\displaystyle{\frac{\partial \tilde{u}_{sh}}{\partial\nu}=\psi\quad\mbox{on\ \ $\partial\Omega$}},\\
&\displaystyle{\frac{\partial \tilde{u}_{sh}}{\partial\nu}=0\quad\mbox{on\ \ $\partial D_\rho$}}.
\end{cases}
\end{equation}
Moreover, we know
\[
u_{sh}=\tilde{u}_{sh}\quad\mbox{and}\quad u=\tilde{u}\ \ \mbox{on\ \ $\partial\Omega$}.
\]

Let $u_0$ be the solution of (\ref{eq:free space}). By straightforward verification, we first see that $Q=u_0-\tilde{u}_{sh}\in H^1(\Omega\backslash\overline{D}_\rho)$ satisfies
\begin{equation}\label{eq:c2}
\begin{cases}
& \displaystyle{\Delta Q+\omega^2 Q=0\quad\mbox{in\ \ $\Omega\backslash\overline{D}_\rho$}},\\
& \displaystyle{\frac{\partial Q}{\partial \nu}=0\quad\mbox{on\ \ $\partial\Omega$}},\\
& \displaystyle{\frac{\partial Q}{\partial \nu}=\frac{\partial u_0}{\partial\nu}}.
\end{cases}
\end{equation}
By Lemma~\ref{lem:order n}, we have
\begin{equation}\label{eq:ineq1}
\|\tilde{u}_{sh}-u_0\|_{H^{1/2}(\partial\Omega)}=\|Q\|_{H^{1/2}(\partial\Omega)}\leq C\rho^N\|\psi\|_{H^{-1/2}(\partial\Omega)}.
\end{equation}
On the other hand, by Theorem~\ref{thm:1}
\begin{equation}\label{eq:ineq2}
\|\tilde{u}-u_0\|_{H^{1/2}(\partial\Omega)}\leq C\rho^N\|\psi\|_{H^{-1/2}(\partial\Omega)}
\end{equation}
Hence, by (\ref{eq:ineq1}) and (\ref{eq:ineq2}) we have
\[
\begin{split}
&\|u_{sh}-u\|_{H^{1/2}(\partial\Omega)}=\|\tilde{u}_{sh}-\tilde{u}\|_{H^{1/2}(\partial\Omega)}\\
\leq &\|\tilde{u}_{sh}-u_0\|_{H^{1/2}(\partial\Omega)}+\|\tilde{u}-u_0\|_{H^{1/2}(\partial\Omega)}\\
\leq & C\rho^N\|\psi\|_{H^{-1/2}(\partial\Omega)}.
\end{split}
\]
\end{pf}

Theorem~\ref{thm:sh} indicates that for a small $\rho$, the FSH layer with a large density parameter really behaves like a sound-hard layer due to that the exterior wave effects are close to each other. Our near-cloaking scheme by employing such FSH lining is shown to produce significantly enhanced cloaking performances compared to the existing ones in literature. The cloaking construction is assessed within general geometry and arbitrary cloaked contents. The assessment is based on controlling the conormal derivative of the wave field on the exterior boundary of the FSH layer and estimating the exterior boundary effects of sound-hard-like small inclusions. Finally, we would like to remark that our present study could be extended to the near-cloaking of full Maxwell's equations by using the technique developed in this work and the estimates due to small electromagnetic inclusions derived in \cite{AVV}, which will be reported in a future paper. 

\section*{Acknowledgement}

The authors would like to thank the anonymous referee for many constructive comments, which have led to significant improvement on the presentation of the paper. 

\bibliographystyle{model1b-num-names}

\end{document}